\documentclass[11pt,a4,leqno,twoside]{amsart}

\usepackage[english,francais]{babel}
\usepackage{amssymb}
\usepackage{amsmath}
\usepackage{graphicx}
\usepackage{amscd}
\usepackage{url}
\usepackage{vmargin}
\usepackage{subeqnarray}
\usepackage{url}

\newtheorem{theoreme}{Th\'{e}or\`{e}me}[section]

\newtheorem{corollaire}[theoreme]{Corollaire}

\newtheorem{definition}[theoreme]{D\'{e}finition}
\newtheorem*{defref}{D\'efinition}
\newtheorem{definitions}[theoreme]{D\'{e}finitions}
\newtheorem{exemple}[theoreme]{Exemple}

\newtheorem{lemme}[theoreme]{Lemme}
\newtheorem{lemme-def}[theoreme]{Lemme-D\'efinition}

\newtheorem{proposition}[theoreme]{Proposition}
\newtheorem{remarque}[theoreme]{Remarque}
\def\preuve{\smallskip\goodbreak{\it Preuve.~~--~\kern.3em}
     \ignorespaces}%
\def\qedbox{$\square$}%
\def\qed{\ifmmode\qedbox\else\unskip\ \hglue0mm\hfill
     \qedbox\smallskip\goodbreak\fi}%
\def\preuved#1{\smallskip\goodbreak{\it Preuve #1.~~--~\kern.3em}
     \ignorespaces}%
\def\preuvedt#1{\smallskip\goodbreak{\it Preuve du
th\'eor\`eme~\ref{#1}.~~--~\kern.3em}
     \ignorespaces}%
\def\idpreuve{\smallskip\goodbreak{\it Id\'ees de la preuve.~~--~\kern.3em}
     \ignorespaces}%\usepackage[all,dvips]{xy}
\def\idpreuvedt#1{\smallskip\goodbreak{\it Id\'ees de la preuve du
th\'eor\`eme~\ref{#1}.~~--~\kern.3em}
     \ignorespaces}%
\newcommand{\GG}{\mathbb{G}}

\newcommand{\ZZ}{\mathbb{Z}}

\newcommand{\CC}{\mathbb{C}}

\newcommand{\NN}{\mathbb{N}}

\newcommand{\D}{\mathcal{D}}

\newcommand{\F}{\mathcal{F}}

\newcommand{\I}{\mathcal{I}}
\newcommand{\J}{\mathcal{J}}

\renewcommand{\L}{\mathcal{L}}

\renewcommand{\O}{\mathcal{O}}

\newcommand{\T}{\mathcal{T}}

\newcommand{\Y}{\mathcal{Y}}
\newcommand{\Z}{\mathcal{Z}}

\newcommand{\thmref}[1]{\noindent \textbf{Th\'{e}or\`{e}me~\ref{#1}.} }

\title{Le groupo\"ide de Galois de $P_1$ et son irr\'eductibilit\'e}
\author{Guy Casale}
\thanks{This work was done when the author was supported by JSPS Postdoctoral Fellowship for Foreign Researchers (FY2004). The author is 
supported by EIF Marie Curie fellowship.\\
\indent \textbf{Address:} Departament de Matematiques, Edifici C, Universitat Autonoma des Barcelona, 08 193 Bellaterra (Barcelona) Spain \\
 \indent \textbf{E-mail address:} \texttt{casale@picard.ups-tlse.fr} ; \texttt{casale@mat.uab.es}
}

\date{}

\address{Departament de Matematiques \\ 
Edifici C \\ 
Universitat Autonoma de Barcelona \\
08 193 Bellaterra (Barcelona)\\
Spain}
 
\address{ \it URL : \rm \url{http://www.picard.ups-tlse.fr/~casale}} 

\email{casale@picard.ups-tlse.fr}
\email{casale@mat.uab.es}



\begin{document}

\maketitle

\selectlanguage{francais}
\begin{abstract}
{ Dans cet article, nous calculons le groupo\"ide de Galois de la premi\`ere \'equation de Painlev\'e. Nous proposons ensuite une d\'efinition de 
r\'eductibilit\'e pour les feuilletages holomorphes singuliers et montrons que la r\'eductibilit\'e peut se lire sur le groupo\"ide de Galois du feuilletage.
Nous obtenons un r\'esultat d'irr\'eductibilit\'e du feuilletage sous-jacent \`a la premi\`ere \'equation de Painlev\'e.}
\end{abstract}

\selectlanguage{english}
\begin{abstract}
\textbf{(The Galois groupoid of $P_1$ and its irreducibility)}
{In this article, the Galois groupoid of the first Painlev\'e equation is computed. This computation uses \'E. Cartan's classification of structural equations
of pseudogroups acting on $\CC^2$ and the degeneration of the first Painlev\'e equation on an elliptic equation ($y'' = 6y^2$). Moreover a definition of reducibility
for singular holomorphic foliations is proposed and a characterization of reducible foliations on theirs Galois groupoids is given. It is applied to prove the 
foliation-irreducibility of the first Painlev\'e equation.
 }
\end{abstract}

\selectlanguage{francais}

\tableofcontents

\section*{Introduction}

Le groupo\"ide de Galois d'un feuilletage a \'et\'e introduit par B. Malgrange dans \cite{legroupoidedegaloisdunfeuilletage}. Cette d\'efinition
concerne les feuilletages holomorphes (singuliers) sur une vari\'et\'e $\CC$-analytique lisse. Cet objet est 
la ``cl\^oture de Zariski'' du pseudo-groupe d'holonomie du feuilletage dans le sens suivant. Soient $(X, \F)$ une vari\'et\'e $\CC$-analytique lisse 
portant un feuilletage $\F$. Les champs de vecteurs locaux tangents \`a $\F$
forment un faisceau en alg\`ebres de Lie de champs de vecteurs. Le pseudo-groupe de transformations de $X$ engendr\'e 
par ces champs de vecteurs est appel\'e le pseudo-groupe tangent du feuilletage et not\'e $\T an(\F)$. Ce pseudo-groupe intervient notamment
dans la construction des transports holonomes. 
Le pseudo-groupe tangent n'est pas d\'ecrit par un syst\`eme d'\'equations aux d\'eriv\'ees partielles, ceci malgr\`es le fait que le faisceau de 
champs de vecteurs dont il est issu soit d\'ecrit par un syst\`eme d'\'equations aux d\'eriv\'ees partielles lin\'eaires. Ce ph\'enom\`ene est bien connu 
sur les groupes alg\'ebriques o\`u le groupe de Lie int\'egrant une sous-alg\`ebre de Lie n'est pas toujours un groupe alg\'ebrique. La d\'efinition 
propos\'ee par B. Malgrange est bas\'ee sur l'absence de  ``troisi\`eme th\'eor\`eme de Lie'' pour les $\D$-groupo\"ides de Lie 
(pseudo-groupes d\'ecrits par des e.d.p. alg\'ebriques ; d\'efinition \ref{groupoidedelie}). 
La d\'efinition du groupo\"ide de Galois est la suivante (nous renvoyons en \ref{definitiongroupoidedegalois} pour une d\'efinition plus pr\'ecise).\\

\begin{defref}
 Le groupo\"ide de Galois d'un feuilletage est le plus petit $\D$-groupo\"ide de Lie  contenant le pseudo-groupe tangent du feuilletage.
 \end{defref}
\vspace{0.3cm}
Dans \cite{legroupoidedegaloisdunfeuilletage}, B. Malgrange montre que cette d\'efinition g\'en\'eralise le groupe de Galois diff\'erentiel d'une 
connexion lin\'eaire int\'egrable.

Dans le cas d'un germe de feuilletage de codimension un, une \'etude compl\`ete des relations entre le groupo\"ide de Galois, la transcendance des 
int\'egrales premi\`eres et les structures g\'eo\-m\'e\-tri\-ques singuli\`eres transverses au feuilletage est faite dans \cite{codimun}. L'outil principal
de cette \'etude est la construction de suites de Godbillon-Vey m\'eromorphes sp\'eciales.

 \subsection{\'Enonc\'e des r\'esultats}
 
 Dans cet article, nous calculons le groupo\"ide de Galois du feuilletage d\'ecrit par la premi\`ere \'equation de Painlev\'e :
 
 $$
 \frac{d^2 y}{dx^2} = 6 y^2 +x.
 $$ 
 Nous noterons $P_1$ cette \'equation et 
 $$
 X_1 = \frac{\partial }{\partial x} + y'\frac{\partial }{\partial y} + (6y^2 + x)\frac{\partial }{\partial y'}
 $$
 le champ de vecteurs d\'ecrivant le feuilletage associ\'e.
 Nous obtenons le r\'esultat suivant :\\
 
 \thmref{groupoidedegaloisdeP1} 
 %-----------------------------%
 \textit{Le groupo\"ide de Galois de $P_1$ est le groupo\"ide d'invariance de la forme $\gamma  = i_{X_1} dx \wedge dy \wedge dy' $. 
 Ces solutions sont les germes de transformations de $\CC^3$, $\Gamma$, satisfaisant les \'equations $\Gamma^*\gamma = \gamma$.
 }
\vspace{0.3cm}
 
Les feuilletages plus simples que les feuilletages de codimension deux sont les feuilletages donn\'es par des \'equations lin\'eaires et les feuilletages 
de codimension un. Nous dirons qu'un feuilletage de codimension deux est r\'eductible si on peut construire deux int\'egrales premi\`eres locales
en utilisant seulement des int\'egrales premi\`eres locales de feuilletages plus simples (d\'efinition \ref{reductibilite}). Plus pr\'ecisement, pour la premi\`ere \'equation de Painlev\'e, 
le th\'eor\`eme \ref{groupoidedegaloisdeP1} nous permet d'obtenir le r\'esultat d'irr\'eductibilit\'e suivant :\\
 
 \thmref{irreductibilite}
 %-------------------%
 \textit{Il n'existe pas d'extension diff\'erentielle $K_n$ de $\left( \CC(x,y,y'), \frac{\partial}{\partial x}, \frac{\partial}{\partial y}, \frac{\partial}{\partial y'} \right)$ 
 contenant une int\'egrale premi\`ere de $X_1$  construite de la mani\`ere suivante :
 $$
 \CC(x,y,y') = K_0 \subset K_1 \ldots \subset K_n 
 $$
 avec 
 \begin{itemize}
     \item $K_{i+1}$ alg\'ebrique sur $K_i$,
     \item ou $K_{i+1} = K_i(h_1,\ldots,h_p)$ avec $dh_j = \sum h_k \omega^k_j$, $\omega^k_j \in K_i \otimes \Omega^1_{\CC^3}$ 
                                                                    et $d\omega^k_j= - \sum \omega^k_\ell \wedge \omega^\ell_j$,  
     \item ou $K_{i+1} = K_i(<h>)$ avec $dh \wedge \omega=0$, $\omega \in K_i \otimes \Omega^1_{\CC^3}$ et $\omega \wedge d\omega =0$.
 \end{itemize}
  }
\vspace{0.3cm}
Le corps $K(h)$ d\'esigne le corps engendr\'e par $h$ et $K(<h>)$ le corps diff\'erentiel engendr\'e par $h$.
Le th\'eor\`ereme que nous montrons est en fait un peu plus g\'en\'eral que celui-ci mais n\'ecessite la mise en place du vocabulaire appropri\'e.

\subsection{Th\'eories de Galois diff\'erentielles et irr\'eductibilit\'e}

\`A la fin du dix-neuvi\`eme si\`ecle, les id\'ees d'\'E. Galois ont \'et\'e \'etendues aux \'equations diff\'erentielles lin\'eaires par \'E. Picard 
\cite{Picard} puis ont \'et\'e compl\'et\'ees par E. Vessiot \cite{Vessiot-memoire}. Dans les ann\'ees 1950,
E. Kolchin d\'eveloppe cette th\'eorie du point de vue des extensions de corps diff\'erentiels \cite{BB,Kolchin}    

Dans \cite{Drach-these}, J. Drach avance une th\'eorie de Galois pour les \'equations diff\'erentielles non-lin\'eaires. Malgr\'e les erreurs qui invalident la 
plupart de ses d\'efinitions, il donne des indications pour calculer le groupo\"ide de Galois de divers feuilletages 
\cite{Drach-equationdelabalistique,Drach-corpspesant} et notament \cite{Drach-painleve} dont cet article suit la d\'emarche.
Dans \cite{Vessiot-definition}, E. Vessiot esquisse une d\'efinition semblable \`a celle de B. Malgrange. Elle est \`a l'origine de la 
d\'efinition du groupe de Galois infinit\'esimal de H.  \cite{Umemura}. Bien qu'ils aient un ``anc\^etre'' commun, l'\'equivalence du 
groupo\"ide de Galois et du groupe de Galois infinit\'esimal n'est encore qu'une conjecture.

Les premi\`eres tentatives de calcul du groupo\"ide de Galois de la premi\`ere \'equation de Painlev\'e sont d\^us \`a P. Painlev\'e 
\cite{Painleve-note} et J. Drach \cite{Drach-painleve}. Utilisant une d\'efinition bas\'ee sur des r\'esultats \'erron\'es, J. Drach donne les directions 
\`a suivre pour prouver le th\'eor\`eme \ref{groupoidedegaloisdeP1}. En faisant appel \`a la classification locale des 
pseudo-groupes de Lie agissant sur $\CC^2$, \'etablie par S. Lie, il affirme que la nature du groupo\"ide de Galois est d\'ecrit par un 
syst\`eme d'\'equations aux d\'eriv\'ees partielles (dit ``r\'esolvant'') admettant une solution rationnelle. Il affirme ensuite que la d\'eg\'en\'erescence
de $P_1$ sur l'\'equation elliptique $y''=6y^2$ peut servir \`a montrer l'abscence de solutions rationnelles aux \'equations r\'esolvantes. 
Le premier point n'est pas justifi\'e et le deuxi\`eme est entach\'e d'erreurs.\\

Dans les articles pr\'ec\'edemment cit\'es, les auteurs affirment que le th\'eor\`eme \ref{groupoidedegaloisdeP1} implique l'irr\'eductibilit\'e ``absolue''
de l'\'equation sans m\^eme d\'efinir cette irr\'eductibilit\'e. Dans les Le\c{c}ons de Stockholm \cite{leconsdestokholm}, P. Painlev\'e d\'efinit 
une notion de r\'eductibilit\'e d'une solution d'une \'equation diff\'erentielle. Une \'equation est dite r\'eductible si toutes ses solutions le sont. Il 
donne ensuite une caract\'erisation des \'equations sans singularit\'es mobiles du second ordre r\'eductibles : la solution g\'en\'erale d\'epend 
semi-transcendentalement des constantes d'int\'egrations. En d'autres termes, l'\'equation admet une int\'egrale premi\`ere rationnelle en les 
variables d\'ependantes.

Cette d\'efinition est tr\`es restrictive comme le montre P. Painlev\'e dans la remarque 28 de \cite{Painleve-bulletin}. Elle a n\'eanmoins 
l'int\'er\^et de faire appara\^itre les diff\'erences entre la r\'eductibilit\'e d'une \'equation (ou du feuilletage sous-jacent) et celle d'une solution 
particuli\`ere.
L'\'etude de la r\'eductibilit\'e des solutions particuli\`eres des \'equations de Painlev\'e est l'{\oe}uvre de l'\'ecole japonaise. H. Umemura 
\cite{Umemura-painleve} et K. Nishioka \cite{Nishioka} donnent un crit\`ere permettant de trouver les familles \`a un param\`etre de solutions r\'eductibles d'une 
\'equation du second ordre et l'appliquent \`a l'\'etude de la premi\`ere \'equation de Painlev\'e. \`A la suite de ces articles, Murata \cite{murata}, 
Watanabe \cite{watanabe, watanabe2}, Noumi-Okamoto \cite{noumi-okamoto} and Umemura-Watanabe \cite{umemura-watanabe, umemura-watanabe-bis} trouvent les solutions 
r\'eductibles non alg\'ebriques des autres \'equations de Painlev\'e. 

Ces r\'esultats prouvent l'irr\'eductibilit\'e (au sens des Le\c{c}ons de Stokholm) des solutions des \'equations de Painlev\'e pour des valeurs g\'en\'eriques 
des param\`etres. Dans \cite{Painleve-bulletin}, P. Painlev\'e souligne la nature restrictive de sa d\'efinition et pose la question des 
rapports entre une d\'efinition de l'irr\'eductibilit\'e du feuilletage sous-jacent \`a l'\'equation et la tentative de th\'eorie de Galois de J. Drach.\\       

Dans cet article, nous proposons une d\'efinition de feuilletages r\'eductibles. La d\'efinition porte sur les feuilletages de codimension 
deux mais il est facile de l'\'etendre aux feuilletages de codimension quelconque. Elle met l'accent sur la nature de certaines
int\'egrales premi\`eres du feuilletage. Les int\'egrales premi\`eres les moins ``transcendantes'' permettent de comprendre grossi\`erement la 
mani\`ere dont la solution g\'en\'erale d\'epend des constantes d'int\'egration. Un des int\'er\^ets de cette d\'efinition est de pouvoir se lire sur le groupo\"ide de Galois du feuilletage. 
Nous suivons les r\'esultats 
partiels de J. Drach pour calculer le groupo\"ide de Galois de $P_1$ et prouver son irr\'eductibilit\'e.    

Dans l'\'etat actuel, la r\'eductibilit\'e au sens des feuilletages ne permet pas de d\'efinir la 
r\'eductibilit\'e d'une solution particuli\`ere. D'une mani\`ere plus g\'en\'erale, les relations entre le type de transcendance d'une solution et celui des 
int\'egrales
premi\`eres sont difficiles \`a saisir. Pour les feuilletages de codimension un, un th\'eor\`eme de M. Singer \cite{singer} dit que si une feuilletage du plan
admet une solution Liouvillienne alors soit elle est alg\'ebrique soit le feuilletage admet une int\'egrale premi\`ere Liouvillienne.

\subsection{Organisation de l'article}

Dans la premi\`ere partie, nous rappelons les d\'efinitions n\'ecessaires \`a la compr\'ehension des r\'esultats et des preuves pr\'esent\'es.  
Nous commen\c{c}ons par un rapide rappel de la construction des espaces de jets et des plus importantes des propri\'et\'es de leurs 
sous-vari\'et\'es. Pour plus de d\'etails, nous renvoyons
aux livres de J.F. Ritt \cite{differentialalgebra} et de J-F. Pommaret \cite{differentialgaloistheory} dans le cadre alg\'ebrique
et \`a B. Malgrange \cite{systemesdifferentiels} dans le cadre analytique.

Nous rappelons ensuite les d\'efinitions de $\D$-groupo\"ides de Lie et de leurs alg\`ebres de Lie. Ces objets formalisent les notions de pseudo-groupes 
de Lie alg\'ebriques de transformations et de leurs alg\`ebres de Lie de champs de vecteurs. Ils sont \'etudi\'es depuis S. Lie et
\'E. Cartan sous le nom de ``groupes infinis de Lie'' ou de ``pseudo-groupes de Lie'' sous des hypoth\`eses suppl\'ementaires de r\'egularit\'e.
L'\'etude de la structure d'un $\D$-groupo\"ide de Lie passe par la compr\'ehension de sa forme de Maurer-Cartan. Cette forme g\'en\'eralise 
la forme de Maurer-Cartan d'un groupe de Lie et satisfait des \'equations de structures semblable \cite{Guillemin-Sternberg, s-s}. 

Enfin, nous donnons la d\'efinition du groupo\"ide de Galois d'un feuilletage.  Nous d\'efi\-nis\-sons les formes invariantes transverses 
du groupo\"ide de Galois d'un feuilletage et montrons comment construire des suites de Godbillon-Vey ``sp\'eciales'' \`a partir de ces formes.
  Cette partie se termine avec l'\'enonc\'e d'un r\'esultat (th\'eor\`eme \ref{itchiban}) donnant une description sommaire des suites de Godbillon-Vey possibles pour un feuilletage 
  d\'efini par une 2-forme ferm\'ee. Ce r\'esultat est une cons\'equence de la preuve donn\'e par \'E. Cartan de la classification des 
  pseudo-groupes agissant sur $\CC^2$. Nous redonnons la preuve de Cartan en annexe.\\

Dans la deuxi\`eme partie, nous calculons le groupo\"ide de Galois de $P_1$. D'apr\`es le th\'eor\`eme \ref{itchiban}, le th\'eor\`eme 
\ref{groupoidedegaloisdeP1} est vrai si trois conditions sont v\'erifi\'ees. Premi\`erement, il n'existe pas d'int\'egrale premi\`ere rationnelle 
du feuilletage. Deuxi\`emement, il n'existe pas de feuilletage de codimension un contenant le feuilletage. Troisi\`emement, il n'existe pas de 
$sl_2$-connexion sur le fibre conormal. 
Chacune de ces conditions se r\'e\'ecrit sous la forme d'un syst\`eme d'\'equations aux d\'eriv\'ees partielles ne devant pas avoir de solutions 
alg\'ebriques. Pour montrer que ces \'equations n'ont pas de solutions, nous utilisons la d\'eg\'en\'erescence de $X_1$ sur 
$X_0 = \frac{\partial }{\partial x} + y'\frac{\partial }{\partial y} + 6y^2\frac{\partial }{\partial y'}$ \`a travers la famille $X_\alpha = \frac{\partial }{\partial x} + y'\frac{\partial }{\partial y} + (6y^2 + \alpha^5 x)\frac{\partial }{\partial y'}$.
La famille de champs de vecteurs $X_\alpha$ \'etant triviale pour $\alpha \in \CC - 0$, les \'equations \`a \'etudier se prolongent rationnellement 
au param\`etre. En les d\'eveloppant le long de $\{\alpha = 0\}$,
on obtient une suite de syst\`emes d'\'equations aux d\'eriv\'ees partielles plus simples.
 Chacun des trois syst\`emes d'e.d.p. donn\'ees par le th\'eor\`eme \ref{itchiban} est \'etudi\'e dans un paragraphe propre. \\
 
Dans la derni\`ere partie, nous d\'efinissons une propri\'et\'e de r\'eductibilit\'e pour les feuilletages. Cette propri\'et\'e signifie que l'on peut 
construire un syst\`eme d'int\'egrales premi\`eres du feuilletage en utilisant successivement des int\'egrales premi\`eres de feuilletages plus simples.
Dans le cas d'un feuilletage de codimension deux, les extensions de corps diff\'erentiels que nous consid\`ererons comme plus simples que 
l'extension du corps des fonctions rationnelles par un syst\`eme de deux int\'egrales premi\`eres ind\'ependantes seront :
\begin{itemize}
\item les extensions alg\'ebriques,
\item les extensions fortement normales (au sens de Kolchin \cite{Kolchin}),
\item les extensions par une int\'egrale premi\`ere d'un feuilletage de codimension un,
\item les extension fortement normales relatives en codimension un (voir d\'efinition \ref{modulo} et \cite{cassidy-singer}) 
\end{itemize}  
Ce dernier type d'extension signifie de mani\`ere grossi\`ere que relativement \`a la projection sur $\CC$ donn\'ee par une int\'egrale premi\`ere, 
le feuilletage que l'on consid\`ere est de codimension un et qu'il admet une int\'egrale premi\`ere dans une extension fortement normale. 
Nous utilisons le type d'une extension de corps diff\'erentiels introduit par Kolchin (\cite{Kolchin}) \`a partir d'un analogue diff\'erentiel du polyn\^ome
de Hilbert pour mesurer la taille du groupo\"ide de Galois d'un feuilletage r\'eductible. Nous montrons que les feuilletages r\'eductibles de codimension
 deux ont un groupo\"ide de Galois ``petit'' (\textit{i.e} de type lin\'eaire) et que le  groupo\"ide de Galois de $P_1$ est plus ``gros'' (\textit{i.e.} 
 de type quadratique).     
 
\subsection{Remerciements}

Ce travail a \'et\'e effectu\'e pendant un s\'ejour post-doctoral \`a l'Universit\'e de Tokyo, je remercie particuli\`erement Kazuo Okamoto et Hitedaka Sakai
pour leur accueil. Je remercie \'egalement Hiroshi Umemura pour son invitation \`a Nagoya, son enthousiasme et son int\'er\^et pour ce travail.   
% Ce travail a aussi b\'en\'efici\'e de nombreux commentaires de B. Malgrange.  

%%%%%%%%%%%%%%%%%%%%%%%%%%%%%%%%%%%%%%%%%%%%%%%%%%%%%%%%%%%%%%%%%%%%%%%%%%%%%%%%%%%%%%%%%%%%%%%%%%%%%%%%%%%%%%%%%%%%%%%%%%%%
%%%%%%%%%%%%%%%%%%%%%%%%%%%%%%%%%%%%%%%%%%%%%%%%%%%%%%%%%%%%%%%%%%%%%%%%%%%%%%%%%%%%%%%%%%%%%%%%%%%%%%%%%%%%%%%%%%%%%%%%%%%%
%%%%%%%%%%%%%%%%%%%%%%%%%%%%%%%%%%%%%%%%%%%%%%%%%%%%%%%%%%%%%%%%%%%%%%%%%%%%%%%%%%%%%%%%%%%%%%%%%%%%%%%%%%%%%%%%%%%%%%%%%%%%
%%%%%
%%

%                 SECTION-DEF                                                      SECTION-DEF                                        SECTION-DEF

%***************************************************************% 
\section{Le groupo\"ide de Galois d'un feuilletage alg\'ebrique}
%***************************************************************%

La notion de $\D$-groupo\"ide de Lie introduite par B. Malgrange dans \cite{legroupoidedegaloisdunfeuilletage} formalise l'id\'ee de 
pseudo-groupe de Lie singulier. Ces objets, sous certaines conditions de r\'egularit\'e, ont \'et\'e \'etudi\'es par divers 
auteurs \cite{Kumpera-Spencer, Guillemin-Sternberg-1, Guillemin-Sternberg, s-s} \`a la suite de S. Lie et \'E. Cartan.

Dans cette partie, nous rappelons les r\'esultats \'el\'ementaires sur les $\D$-groupo\"ides de Lie. Pour les d\'etails nous
renvoyons \`a B. Malgrange \cite{legroupoidedegaloisdunfeuilletage} et \`a J-F. Pommaret \cite{differentialgaloistheory}.

%%%%%%%%%%%%%%%%%%%%%%%%%%%%%%%%%%%%%%%%%%%%%%%%%%%%%%%%%%%%%%%%%%%%%%%%%%%%%%%%%%%%%%%%%%%%%%%%%%%%%%%%%%%%%%%%%%%%%%%%%%%%
%%%%%%%%%%%%%%%%%%%%%%%%%%%%%%%%%%%%%%%%%%%%%%%%%%%%%%%%%%%%%%%%%%%%%%%%%%%%%%%%%%%%%%%%%%%%%%%%%%%%%%%%%%%%%%%%%%%%%%%%%%%%
%%

%              SUBSECTION-DEF.1                                              SUBSECTION-DEF.1                                       SUBSECTION-DEF.1

%+++++++++++++++++++++++++++++++++++++%
\subsection{Espaces de jets et $\D$-vari\'et\'es}
%+++++++++++++++++++++++++++++++++++++%

Soit $X$ une vari\'et\'e alg\'ebrique lisse sur $\CC$ et $Z \to X$ une vari\'et\'e alg\'ebrique sur $X$. 
Nous noterons $J_q(Z_{/X})$ l'espace des jets d'ordre $q$ de sections de $Z$ sur $X$. 
Ces espaces sont des espaces affines au-dessus de $Z$. 
Rappelons bri\`evement comment on les construit. 

Commen\c{c}ons par construire $J_q(X \times \CC^N _{/X})$ o\`u $X$ est affine de dimension $n$ telle qu'il existe un 
rev\^etement d'un ouvert de $\CC^n$ par $X$ non ramifi\'e. 
Nous noterons $x_1,\ldots x_n$ des coordonn\'ees sur $\CC^n$ et celles induites sur $X$ et $y_1,\ldots y_N$ celles de $\CC^N$. 
L'espace $J_q(X \times \CC^N _{/X})$ est alors d\'efini comme l'espace $X \times \CC^N$ muni du faisceau d'anneaux 
%
%aaaaaaaaaaa
      $$
      \O_{J_q(X \times \CC^N_{/X})} = \O_{X \times \CC^N}[y_i^\alpha],\ 1 \leq i \leq N,\ \alpha \in \NN^n,\ 1 \leq |\alpha| = \alpha_1 + \ldots + \alpha_n \leq q. 
      $$
%aaaaaaaaaaa
%
Ces espaces de jets sont munis de d\'erivations $D_i$ de $\O_{J_q(X \times \CC^N _{/X})}$ dans $\O_{J_{q+1}(X \times \CC^N _{/X})}$
d\'efinies par 
%
%aaaaaaaaaaa
      $$
      D_j(x_i)= \delta_{ij}\ , \ \ D_j{y_i^\alpha} = y_i^{\alpha + \epsilon_j},
      $$
%aaaaaaaaaaa
%
$\epsilon_j$ \'etant le multi-indice 
de poids 1 dont la seul composante non nulle est la $j$-i\`eme.
  
Soit $Z$ une vari\'et\'e affine au-dessus de $X$ (plong\'ee dans $X \times \CC^N$) d\'ecrite par un id\'eal $I$. 
L'espace des jets $J_q(Z_{/X})$ est le sous-espace de $J_q(X \times \CC^N _{/X})$ d\'efini au-dessus de $Z$ par l'id\'eal
$\sum_{|\alpha| \leq q} D^{\alpha}I$. Nous noterons $\O_{J_q(Z_{/X})}$ le faisceau quotient sur $Z$. 
C'est le faisceau des \'equations aux d\'eriv\'ees partielles d'ordre $q$ portant sur les sections de $Z$ sur $X$. 
Ces espaces ne d\'ependent ni des coordonn\'ees choisies sur $X$ ni du plongement de $Z$. 
Lorsque $X$ et $Z$ ne sont pas affines, on construit les espaces de jets par recollement des constructions locales. 
Ces espaces sont munis de projections naturelles 
%
%aaaaaaaaaa
      $$
      \pi_q^{q+s} : J_{q+s}(Z_{/X}) \to J_q(Z_{/X}).
      $$
%aaaaaaaaaa
Lorsque $Z$ est le fibr\'e trivial $X \times Y$, nous noterons $J_q(X \to Y)$ l'espace des jets de sections de $X \times Y$ sur $X$. 

\begin{remarque}
La construction usuelle des espaces de jets donne des vari\'et\'e alg\'ebriques $\J_q(Z_{/X})$ munit de leurs faisceaux structuraux $\O_{\J_q(Z_{/X})}$.
Ici suivant \cite{legroupoidedegaloisdunfeuilletage}, nous consid\'erons l'espace $Z$ annel\'e par l'image direct du faisceau $\O_{\J_q(Z_{/X})}$.
\end{remarque}
%
%-----------%
\begin{lemme}
%-----------%
\label{jet<->vertical}
%-----------%
   Les constructions des espaces de jets et du tangent relatif (ou ``vertical'') commutent \textit{i.e.} il existe un isomorphisme canonique 
      $$
      T \left(J_q(Z_{/X})\right)_{/X} \simeq J_q\left( (T Z_{/X})_{/X} \right) 
      $$
%----------%   
\end{lemme}
%----------%

\begin{preuve}
%------------%
      Cette formule n'est qu'une version de la commutativit\'e des d\'erivations partielles.
      V\'erifions la lorsque $Z$ est affine au-dessus de $X$. La construction du tangent relatif 
      se fait de la mani\`ere suivante. On plonge $Z$ dans $X \times \CC^N$ et on note $I$ l'id\'eal d\'efinissant l'image 
      de $Z$. Soient $z_1,\ldots,z_N$ des coordonn\'ees sur $\CC^N$, le tangent relatif de $X \times \CC^N$ est l'espace vectoriel
      trivial sur cet espace : $(X \times \CC^N) \times \CC^N$ avec les coordonn\'ees $dz_1, \ldots, dz_N$ sur les fibres. 
      Pour $f \in \O_X[z_1,\ldots,z_N]$ on notera $ df$ la fonction $\sum \frac{\partial f}{\partial z_i} dz_i$. 
      Le tangent relatif de $Z$ sur $X$ est l'espace vectoriel sur $Z$ d\'efini par l'id\'eal $I + dI$.
     
       Lorsque $X$ est un rev\^etement non ramifi\'e d'un ouvert de $\CC^n$, l'espace $J_q((TZ_{/X})/X)$ est d\'efini par les \'equations suivantes :
      $$
      D_j d \, f = \sum_i \frac{\partial^2 f}{\partial x_j \partial z_i} dz_i + \sum_{i,k} \frac{\partial^2 f}{\partial z_k \partial z_i} z_k^{\epsilon_j}dz_i + \sum_i \frac{\partial f}{\partial z_i} (dz_i)^{\epsilon_j} 
      $$
      pour toutes les \'equations $f \in I$. D'un autre cot\'e, l'espace $T(J_q(Z_{/X}))_{/X}$ est d\'efinie par les \'equations :
      $$
      d \, D_j f =\sum_i \frac{\partial^2 f}{\partial x_j \partial z_i} dz_i + \sum_{i,k} \frac{\partial^2 f}{\partial z_k \partial z_i} z_k^{\epsilon_j}dz_i + \sum_i \frac{\partial f}{\partial z_i} d(z_i^{\epsilon_j}).  
      $$ 
      L'isomorphisme pour $q=1$ est donn\'e par $(dz_i)^{\epsilon_j} = d(z_i^{\epsilon_j})$. Des formules analogues donnent les isomorphismes pour $q>1$.     
%-------------%
\end{preuve}

%
%----------------%
\begin{definition}
%----------------%
\label{prolongement}
%----------------%
   Soit $Y_q$ une sous-vari\'et\'e de $J_q(Z_{/X})$ donn\'ee par un id\'eal $I_q$. Le prolongement d'ordre 1 de $Y_q$ 
   est la sous-vari\'et\'e de $J_{q+1}(Z_{/X})$ d\'efinie par l'id\'eal $I_q + \sum D_i I_q$. Cette vari\'et\'e sera 
   not\'ee $pr_1 Y_q$ et son id\'eal $pr_1 I_q$. On d\'efinit par it\'eration les prolongements d'ordre $k$, $pr_k I_q$.
%---------------%
\end{definition}
%---------------%    
%
%
%---------------%
\begin{remarque}
   On n'a pas forcement $\pi_q^{q+s}pr_s I_q \not = I_{q}$.
\end{remarque}
%---------------%

On d\'efinit l'espace des jets (d'ordre infini) de sections de $Z$ sur $X$ comme la pro-vari\'et\'e affine sur $Z$
%
%aaaaaaaaaa
      $$
      J(Z_{/X}) = \underset{\leftarrow}{\lim}\ J_q(Z_{/X}),
      $$
%aaaaaaaaaa
%
c'est-\`a-dire  
$Z$ annel\'ee par $\O_{J(Z_{/X})} = \underset{\rightarrow}{\lim}\ \O_{J_q(Z_{/X})}$. Ces anneaux sont filtr\'es par l'ordre 
des \'equations et muni d'une connexion 
%
%aaaaaaaaaa
      $$
      D :  \O_{J(Z_{/X})} \to \O_{J(Z_{/X})}\otimes_{\O_X}\Omega^1_X  
      $$ 
%aaaaaaaaaa
%
d\'efinie en coordonn\'ees locales par $D f = \sum  D_i f \otimes dx_i $.
 
On a alors la d\'efinition suivante de ``\textit{syst\`emes d'\'equations diff\'erentielles portant sur les sections de $Z$ sur $X$}''
appel\'es par la suite $\D$-sous-vari\'et\'es de $J(Z_{/X})$ ou plus bri\`evement $\D$-vari\'et\'es.
%
%----------------%
\begin{definition}
%----------------%
\label{dvariete}
%----------------%
   Une $\D$-(sous)-vari\'et\'e de $J(Z_{/X})$ est une sous-vari\'et\'e $\Y$ d\'efinie par un 
   faisceau d'id\'eaux $\I$ de $\O_{J(Z_{/X})}$ tel que :
   \begin{enumerate}
      \item $\I$ est \textit{pseudo-coh\'erent}, \textit{i.e.} les id\'eaux $\I_q = \I \cap \O_{J_q(Z_{/X})}$ sont coh\'erents,
      \item $\I$ est diff\'erentiel, \textit{i.e.} $D \I \subset \I \Omega^1_X$.
   \end{enumerate}
%----------------%   
\end{definition}
%----------------%
%
Le faisceau d'anneaux de $\Y$, $\O_{\Y}= \O_{J^*(Z_{/X})}/\I$ est donc muni d'une connexion induite 
$D : \O_\Y \to \O_\Y \otimes_{\O_X}\Omega^1_X$.
On pourra d\'efinir les $\D$-vari\'et\'es g\'en\'erales : on recollera des espaces d\'efinis comme ci-dessus en respectant 
les connexions. Nous n'en aurons pas besoin. 

Une solution (resp. formelle) de $\Y$ en un point $x \in X$ est un morphisme $f$ de $\O_\Y$ dans $\O^{an}_{X,x}$ 
(resp. $\widehat{\O}_{X,x}$) au-dessus de la restriction $\O_X \to \O^{an}_{X,x}$ (resp. $\widehat{\O}_{X,x}$)
et diff\'erentiel \textit{i.e.} $d \circ f = f \circ D$. Les solutions correspondent aux sections de $Z$ 
v\'erifiant les \'equations diff\'erentielles d\'efinissant $\Y$.

%Les solutions formelles s'identifient au points de $\Y$.

Pour finir, rappelons les r\'esultats suivants sur les $\D$-vari\'et\'es. 
%
%-----------%
\begin{lemme}[Ritt \cite{differentialalgebra}]
%-----------%
   Si $\I$ est un id\'eal diff\'erentiel alors son radical $\I^{red}$ est aussi diff\'erentiel.
%-----------%
\end{lemme}
%-----------%
%
%
%-----------%
\begin{lemme}[\cite{legroupoidedegaloisdunfeuilletage}]
%-----------%
   Si $\Y$ est une $\D$-vari\'et\'e r\'eduite, $\O_{\Y}$ est sans torsion sur $\O_X$.
%-----------%
\end{lemme}
%-----------%
%
%
%--------------%
\begin{theoreme}[Ritt-Raudenbush \cite{differentialalgebra}]
%--------------%
   Soit $\Y$ une $\D$-vari\'et\'e r\'eduite d\'efinie par $\I$ et $\I_q = \I \cap \O_{J^*_q(Z_{/X})}$ l'id\'eal des 
   \'equations d'ordre $q$.
   Il existe un entier $q$ tel que $\I$ soit l'id\'eal diff\'erentiel r\'eduit engendr\'e par $\I_q$.
%--------------%
\end{theoreme}
%--------------%
%
%
%--------------%
\begin{theoreme}[\cite{differentialalgebra}]
%--------------%
   Soit $\Y$ une $\D$-vari\'et\'e r\'eduite d\'efinie par $\I$ et $\Y_q$ l'espace d\'efini par $\I_q = \I \cap \O_{J^*_q(Z_{/X})}$.
   Il existe un entier $q$ et une hypersurface $S \subset \Y_q$ tels que par chaque points de $\Y_q-S$ 
   passe une solution convergente.
%-------------% 
\end{theoreme}
%-------------%
%
Ces th\'eor\`emes ont \'et\'e g\'en\'eralis\'es par B. Malgrange dans la situation mixte d'une vari\'et\'e analytique $Z$ 
sur $X$ analytique lisse et $\Y$ d\'efinie par des \'equations aux d\'eriv\'ees partielles polynomiales en les 
d\'eriv\'ees. Dans cette situation, B. Malgrange montre un r\'esultat \textit{a priori} plus fort que celui de Ritt que
 nous rappelons ici dans le cas ``tout alg\'ebrique''.
%
%--------------% 
\begin{theoreme}[d'involutivit\'e g\'en\'erique \cite{differentialgaloistheory, systemesdifferentiels}]   
%--------------%
   Soit $\Y$ une $\D$-vari\'et\'e r\'eduite d\'efinie par $\I$ et $\Y_q$ l'espace d\'efini par $\I_q = \I \cap \O_{J^*_q(Z_{/X})}$.
   Il existe un entier $q$ et une hypersurface $S \subset \Y_q$ tels que $\Y_q$ soit involutif en dehors de $S$.
%-------------% 
\end{theoreme}
%-------------%
%
L'involutivit\'e est une condition de r\'egularit\'e d'un syst\`eme d'\'equation aux d\'eriv\'ees partielles sous laquelle on a 
un th\'eor\`eme d'existence de solutions convergentes : le th\'eor\`eme de Cartan-K\"ahler. Nous ne rappelerons pas cette 
condition et renvoyons aux r\'ef\'erences du th\'eor\`eme.

%%%%%%%%%%%%%%%%%%%%%%%%%%%%%%%%%%%%%%%%%%%%%%%%%%%%%%%%%%%%%%%%%%%%%%%%%%%%%%%%%%%%%%%%%%%%%%%%%%%%%%%%%%%%%%%%%%%%%%%%%%%%
%%%%%%%%%%%%%%%%%%%%%%%%%%%%%%%%%%%%%%%%%%%%%%%%%%%%%%%%%%%%%%%%%%%%%%%%%%%%%%%%%%%%%%%%%%%%%%%%%%%%%%%%%%%%%%%%%%%%%%%%%%%%
%%

%          SUBSECTION-DEF.2                                      SUBSECTION-DEF.2                                       SUBSECTION-DEF.2

%++++++++++++++++++++++++++++++++++++++++++++++++++++++%
\subsection{$\D$-groupo\"ides de Lie alg\'ebriques et $\D$-alg\`ebres de Lie}
%++++++++++++++++++++++++++++++++++++++++++++++++++++++%

\subsubsection{$\D$-groupo\"ides de Lie alg\'ebriques}
%------------------------------------------------%
Nous allons maintenant rappeler la d\'efinition de $\D$-groupo\"ide de Lie alg\'ebrique sur une vari\'et\'e alg\'ebrique 
lisse $X$ sur $\CC$. L'espace $J_q(X \to X)$  et son ouvert $J^*_q(X \to X)$ des jets inversibles (d\'efini par $det(y_i^{\epsilon_j}) \not = 0$) 
seront not\'es $J_q$ et $J^*_q$ lorsqu'aucune confusion ne pourra \^etre faite sur $X$. 

Ces espaces sont naturellement munis de deux projections sur $X$, la \textit{source} $s$ et le \textit{but} $t$, 
correspondant aux projections de $X \times X$ sur le premier et le deuxi\`eme facteur, d'une \textit{composition} 
partielle associative 
%aaaaaaaaa
      $$
      c : \left(J_q, t\right) \times_X \left(J_q, s\right) \to J_q
      $$
%aaaaaaaaa
induite par la composition des applications formelles de $X$ dans elle-m\^eme et d'une application \textit{identit\'e}
%aaaaaaaaa
      $$
      e : X \to J_q
      $$
%aaaaaaaaa
donn\'ee par les jets d'ordre $q$ de l'identit\'e. L'espace des jets inversibles est de plus muni d'une \textit{inversion} 
%aaaaaaaaa
      $$
      i : J^*_q \to J^*_q
      $$
%aaaaaaaaa
induite par l'inversion des diff\'eomorphismes formels de $X$. Ces applications font de $J^*_q$ un groupo\"ide 
agissant sur $X$.
Elles sont compatibles aux projections naturelles $\pi_q^{q+s}$
et donne une structure de groupo\"ide sur $J^* = \underset{\leftarrow}{\lim}\ J_q$, c'est-\`a-dire sur 
$X \times X$ muni de l'anneau $\O_{J^*} = \underset{\rightarrow}{\lim}\ \O_{J^*_q}$. 
%
%----------------% 
\begin{definition}
%----------------%
\label{groupoide}
%----------------%
   Un sous-groupo\"ide de $J^*_q$ est une sous-vari\'et\'e $Y_q$ donn\'ee par un faisceau coh\'erent d'id\'eaux $I_q$ 
   satisfaisant :
   \begin{enumerate}
      \item $I_q \subset \ker e^*$,
      \item $i^*I_q \subset I_q$,
      \item $c^*I_q \subset \pi_{1}^*I_q + \pi_{2}^*I_q$, 
   \end{enumerate}  
   $\pi$ d\'esignant les projections de $\left(J_q, t\right) \times_X \left(J_q, s\right)$ sur $J_q$.
%--------------%
\end{definition}
%--------------%
% 
On a la version alg\'ebrique d'un r\'esultat classique de g\'eom\'etrie diff\'erentielle.
%
%-----------------%
\begin{proposition}[\cite{legroupoidedegaloisdunfeuilletage}]
%-----------------%
   Si $Y_q$ un sous-groupo\"ide de $J^*_q$ alors $pr_1 Y_q$ est un sous-groupo\"ide de $J^*_{q+1}$.
%-----------------%
\end{proposition}
%-----------------% 
%
%
%--------------%
\begin{remarque}
   Soient $Y_q$ un sous-groupo\"ide d\'efinie par $I_q$ et $\Y$ la $\D$-vari\'et\'e d\'efinie par 
   $\I = \cup pr_s I_q$. Notons $\Y_q$ la sous-vari\'et\'e d\'efinie par $\I \cap \O_{J^*_q}$. \textit{A priori}
    $\Y_q$ n'est pas un sous-groupo\"ide de $J^*_q(X \to X)$ mais seulement en dehors d'une sous-vari\'et\'e $S$ 
    \textit{i.e.} $\Y_q$ est un sous-groupo\"ide de $J^*_q(X-S \to X-S)$
\end{remarque} 
%--------------%
 Cette remarque montre que la d\'efinition naturelle suivante :\textit{ ``une $\D$-vari\'et\'e $\Y$ est un sous-groupo\"ide
 de $J^*$ si il existe $k$ tel que pour tout $q \geq k$ $\Y_q$ est un sous-groupo\"ide de $J^*_q$''} est tr\`es restrictive. De plus certains 
 pseudo-groupes de transformations, comme les pseudo-groupes d'invariances de fonctions rationnelles, ne rentrent dans le cadre 
 de cette d\'efinition qu'en dehors d'une sous-vari\'et\'e ferm\'ee de codimension au moins un (dans le cas pr\'ec\'edent : 
 le lieu d'ind\'etermination).

%
%----------------% 
\begin{definition}
%----------------%
\label{groupoidesing}
%----------------%
   Soit $S$ une hypersurface de $X$. Un sous-groupo\"ide de $J^*_q(X \to X)$ \`a singularit\'es sur $S$
   est une sous-vari\'et\'e ferm\'ee donn\'ee par un faisceau coh\'erent d'id\'eaux $I_q$ telle que
   la trace de la vari\'et\'e sur $J^*_q(X-S \to X-S)$ d\'ecrive un sous-groupo\"ide.  
%--------------%
\end{definition}
%--------------%
% 
 
%
%----------------% 
\begin{definition}
%----------------%
\label{groupoidedelie}
%----------------%
   Un $\D$-groupo\"ide de Lie sur $X$ est une $\D$-vari\'et\'e r\'eduite $\Y$ de $J^*$donn\'ee par un id\'eal $\I$ v\'erifiant :
   il existe une hypersurface $S$ de $X$ et un entier $k$ tels que
   pour tout $q \geq k$ , $\Y_q$ est un sous-groupo\"ide de $J^*_q$ \`a singularit\'es sur $S$. 
%--------------%
\end{definition}
%--------------%
% 
 
%--------------%
\begin{remarque}
   Dans les conditions de la d\'efinition pr\'ecedente, pour tout entier $q$, 
   $\Y_q$ est un sous-groupo\"ide singulier. En effet fixons un entier $q_0$ tel que $\Y_{q_0}$ soit un groupo\"ide singulier 
   sur $S$ et prenons $q < q_0$. Par d\'efinition $\pi_{q}^{q_0}$ est dominant donc surjectif en dehors d'un ensemble 
   alg\'ebrique $S'$ de $\Y_q$, on a donc 
      $$
      c \left( \Y_q-(S'\cup \pi_1^{-1} S \cup \pi_2^{-1} S) \times_{X} \Y_q-(S'\cup \pi_1^{-1} S \cup \pi_2^{-1} S) \right)
       \subset \pi_q^{q_0}\Y_{q_0} \subset \Y_q.
      $$
   Des arguments classiques (voir par exemple la preuve du lemme $4.3.3$ de \cite{legroupoidedegaloisdunfeuilletage}) montrent alors que $\Y_q$ 
   est un groupo\"ide \`a singularit\'es sur une vari\'et\'e $S''$ contenant $S$.
\end{remarque}
%--------------%

%%%%%%%%%%%%%%%%%%%%%%%%%%%%%%%%%%%%%%%%%%%%%%%%%%%%%%%%%%%%%%%%%%%%%%%%%%%%%%%%%%%%%%%%%%%%%%%%%%%%%%%%%%%%%%%%%%%%%%%%%%%%
%%%%%%%%%%%%%%%%%%%%%%%%%%%%%%%%%%%%%%%%%%%%%%%%%%%%%%%%%%%%%%%%%%%%%%%%%%%%%%%%%%%%%%%%%%%%%%%%%%%%%%%%%%%%%%%%%%%%%%%%%%%%
%%

%\subsection{Invariants diff\'erentielles d'un $\D$-groupo\"ide de Lie alg\'ebrique}
%---------------------------------------------------------------------------------%

Nous allons rappeler quelques r\'esultats sur l'action d'un $\D$-groupo\"ide de Lie $\Y$ sur $J^*$ par composition au but.
Nous noterons $c_\Y$ la restriction de $c$ sur $\left(J^* , t \right) \times_{X} \left( \Y , s \right)$ \`a valeurs dans $J^*$.  

%
%----------------%
\begin{definition}
%----------------%
\label{invariantsdifferentiels}
%----------------%
   Un invariant diff\'erentiel d'ordre $q$ pour $\Y$ est une fonction rationnelle $F$ sur $J^*_q$ telle que 
   $F \circ c_\Y = F \circ \pi$ ($\pi$ d\'esignant la premi\`ere projection de $J^*_q \times_{X} \Y_q$ sur $J^*_q$).
%----------------%
\end{definition}
%----------------%
%
\begin{remarque}
   N'importe quelle fonction sur $X$ que l'on remonte sur $J^*_q$ par la projection source est un invariant.
\end{remarque}
On dira qu'un ensemble d'invariants diff\'erentiels $\{F_i\}_{1 \leq i \leq p}$ d\'efinis sur un ouvert de Zariski $U$ forment un syst\`eme 
complet d'invariants pour $\Y$ si $\Y$ est la $\D$-vari\'et\'e r\'eduite d\'efinie par l'adh\'erence de Zariski de la $\D$-vari\'et\'e d\'efinie par 
les \'equations $F_i - F_i \circ e \circ s$ au-dessus de $U$.

Le th\'eor\`eme suivant est prouv\'e d'une mani\`ere analogue au th\'eor\`eme de Chevalley-Kolchin dans (\cite{differentialgaloistheory}, 467--469
  proposition 2.36 du chapitre 3). Il r\'esulte aussi de l'existence d'un quotient g\'en\'erique pour une relation d'\'equivalence (\cite{SGA}, th\'eor\`eme 8.1).
%--------------%
\begin{theoreme}[\cite{differentialgaloistheory, SGA}]
%--------------%
\label{theoinvdiff}
%--------------%
   Si $\Y$ un $\D$-groupo\"ide de Lie agissant sur un vari\'et\'e $X$ alors il admet un syst\`eme complet 
   d'invariants diff\'erentiels.
%--------------%
\end{theoreme}
%--------------%
Dans le contexte de cet article, le th\'eor\`eme 8.1 de \cite{SGA} s'applique comme suit. On se place sur $J^*_q$ ($q$ assez grand pour 
que $\Y_q$ satisfasse aux hypoth\`eses du th\'eor\`eme de Ritt-Raudenbush) et on consid\`ere la relation d'\'equivalence  $\alpha \sim \beta$ si $s(\alpha) = s(\beta)$
et $\beta \circ \alpha^{-1} \in \Y_q$. Le r\'esultat pr\'ec\'edent donne l'existence d'un ouvert de Zariski $U \subset J^*$ et d'une vari\'et\'e quotient de $U$ par
la relation d'\'equivalence. En tirant en arri\`ere des fonctions rationnelles du quotient, on obtient un syst\`eme complet d'invariants diff\'erentiels. 
%\begin{preuve}

%\end{preuve}

%%%%%%%%%%%%%%%%%%%%%%%%%%%%%%%%%%%%%%%%%%%%%%%%%%%%%%%%%%%%%%%%%%%%%%%%%%%%%%%%%%%%%%%%%%%%%%%%%%%%%%%%%%%%%%%%%%%%%%%%%%%%
%%%%%%%%%%%%%%%%%%%%%%%%%%%%%%%%%%%%%%%%%%%%%%%%%%%%%%%%%%%%%%%%%%%%%%%%%%%%%%%%%%%%%%%%%%%%%%%%%%%%%%%%%%%%%%%%%%%%%%%%%%%%
%%
%%
\subsubsection{$\D$-alg\`ebres de Lie}
%-----------------------------------%
\label{repere}

Une $\D$-alg\`ebre de Lie sur $X$ est un faisceau d'alg\`ebres de Lie de champs de vecteurs d\'efini par un syst\`eme 
d'\'equations aux d\'eriv\'ees partielles. La d\'efinition que nous allons rappeler met l'accent sur les \'equations 
d\'efinissant l'alg\`ebre.

Un espace vectoriel $V$ sur $X$ est d\'efini de la mani\`ere suivante. Lorsque $X$ est affine, un espace vectoriel est un sous-espace d'un 
espace vectoriel trivial $X \times \CC^n$ d\'efini par des \'equations lin\'eaires sur les fibres. Lorsque $X$ n'est pas affine, on demande que $V$
soit localement de la forme pr\'ec\'edente et que les recollements soient lin\'eaires sur les fibres. Cette d\'efinition autorise les fibres \`a changer
de dimension.  
  
Les espaces $J_q(TX_{/X})$ sont des espaces vectoriels sur $X$. Nous allons construire un crochet sur les sections de ces 
espaces g\'en\'eralisant le crochet de Lie des sections de $TX$. Ce crochet est le crochet de Spencer et la construction 
que nous allons donner est celle de \cite{legroupoidedegaloisdunfeuilletage}.

Nous noterons $R_q(X)$ (ou juste $R_q$) l'espace des rep\`eres d'ordre $q$ sur $X$ qui s'identifie au sous-espace de $J^*_q$
d\'efini en fixant une valeur \`a l'application source. Cet espace est un fibr\'e principal sur $X$ de groupe structural 
 le groupe des jets d'ordre $q$ de biholomorphismes de $(\CC^n,0)$ not\'e $\Gamma^n_q$. L'application 
$\lambda : R_q \times R_q \to J^*_q$ qui \`a deux rep\`eres $r_q$ et $s_q$ associe le jet $ s_q \circ r_q^{-1}$ s'identifie 
au quotient sous l'action diagonale de $\Gamma^n_q$ sur les deux facteurs. Cette application induit une application du tangent relatif
$T(R_q \times R_q)_{/R_q}|_{diag}$ le long de la diagonale sur le tangent vertical $T(J^*_q)_{/X}|_{id}$ 
le long de l'identit\'e. Le lemme \ref{jet<->vertical} identifie ce dernier espace vectoriel \`a $J_q(TX_{/X})$ et permet 
de projeter le crochet de Lie de $T R_q$ sur un crochet sur les sections de $J_q(TX_{/X})$. C'est le crochet de Spencer. 
%
%----------------%
\begin{definitions}
%----------------%
   Une sous-alg\`ebre de Lie de $J_q(TX_{/X})$ est un sous-espace vectoriel dont les sections analytiques locales sont stables sous le crochet
   de Spencer.
   
   Un $\D$-(sous)-espace vectoriel $\L$ de $J(TX_{/X})$ est une $\D$-vari\'et\'e d\'efinie par des \'equations lin\'eaires.
   
   Une $\D$-(sous)-alg\`ebre de Lie $\L$ de $J(TX_{/X})$ est un $\D$-espace vectoriel tel que les espaces $\L_q$ soient
   des sous-alg\`ebres de Lie de $J_q(TX_{/X})$.
  
%---------------% 
\end{definitions}
%---------------%
%
On construit une alg\`ebre de Lie $LY_q$ \`a partir d'un groupo\"ide de Lie singulier $Y_q$ en consid\'erant le tangent 
relatif \`a la projection source ${TY_q}_{/X}|_{id}$ ou son image dans $J_q(TX_{/X})$ (par le lemme \ref{jet<->vertical}). 
La stabilit\'e des sections du lin\'earis\'e sous le crochet de Spencer est une cons\'equence de la stabilit\'e de la vari\'et\'e $Y_q$ 
sous la composition. Nous renvoyons \`a \cite{legroupoidedegaloisdunfeuilletage} pour les d\'etails ainsi que pour la preuve de la 
proposition suivante. 

%
%-----------------%
\begin{proposition}[\cite{legroupoidedegaloisdunfeuilletage}]
%-----------------%
Si $\Y$ est un $\D$-groupo\"ide de Lie, son tangent vertical le long de l'identit\'e ${T\Y}_{/X}|_{id} \subset J(TX_{/X})$
est une $\D$-alg\`ebre de Lie. Elle sera not\'ee $L\Y$.
%------------------%
\end{proposition}
%------------------%
%

\subsubsection{Le groupo\"ide de Galois d'un feuilletage}

\begin{definition}   
   Un feuilletage (singulier) de $X$ est une $\D$-alg\`ebre de Lie $\F$ sans torsion sur $X$ et diff\'erentiellement engendr\'e par $\F_0$.
\end{definition}
    
%
%--------------%
\begin{remarque}
%--------------%
\label{remarquenonintegrable}
   Une $\D$-alg\`ebre de Lie n'est pas toujours la $\D$-alg\`ebre de Lie d'un $\D$-groupo\"ide de Lie. Prenons par exemple un 
   feuilletage $\F$ de codimension $d$. D'apr\`es le th\'eor\`eme \ref{theoinvdiff}, si il existe un $\D$-groupo\"ide de Lie 
   dont $\F$ est la $\D$-alg\`ebre de Lie alors le feuilletage admet $d$ int\'egrales premi\`eres rationnelles ind\'ependantes. 

   Un feuilletage sans int\'egrales premi\`eres rationnelles est donc un exemple de $\D$-alg\`ebre de Lie ne provenant pas d'un
   $\D$-groupo\"ide de Lie.
%--------------%
\end{remarque}    
%--------------%
%
    
La d\'efinition suivante  g\'en\'eralise le groupe de Galois diff\'erentiel d'une \'equation lin\'eaire. 
%
%----------------%
\begin{definition}[\cite{legroupoidedegaloisdunfeuilletage}]
%----------------%
\label{definitiongroupoidedegalois}
   Soit $\F$ un feuilletage sur $X$. Le groupo\"ide de Galois de $\F$ est le plus petit $\D$-groupo\"ide de Lie dont la 
   $\D$-alg\`ebre de Lie contient le feuilletage. Il sera not\'e $Gal(\F)$.
%----------------%
\end{definition}
%----------------%
%
Autrement dit $Gal(\F)$ est le plus petit $\D$-groupo\"ide de Lie contenant le pseudo-groupe $\T an (F)$. Par commodit\'e, nous posons la d\'efinition suivante. 
%
%----------------%
\begin{definition}
%----------------%
   Soit $\F$ un feuilletage sur $X$. Un $\D$-groupo\"ide de Lie pr\'eservant le feuilletage dont la $\D$-alg\`ebre de Lie contient le feuilletage 
   sera dit admissible pour $\F$.
%----------------%
\end{definition}
%----------------%
%    

%%%%%%%%%%%%%%%%%%%%%%%%%%%%%%%%%%%%%%%%%%%%%%%%%%%%%%%%%%%%%%%%%%%%%%%%%%%%%%%%%%%%%%%%%%%%%%%%%%%%%%%%%%%%%%%%%%%%%%%%%%%%
%%%%%%%%%%%%%%%%%%%%%%%%%%%%%%%%%%%%%%%%%%%%%%%%%%%%%%%%%%%%%%%%%%%%%%%%%%%%%%%%%%%%%%%%%%%%%%%%%%%%%%%%%%%%%%%%%%%%%%%%%%%%
%%

%           SUBSECTION-DEF.3                                    SUBSECTION-DEF.3                                      SUBSECTION-DEF.3

%++++++++++++++++++++++++++++++++++++++++++++++++%
\subsection{Forme de Maurer-Cartan d'un $\D$-groupo\"ide de Lie et suites de Godbillon-Vey d'un feuilletage}
%++++++++++++++++++++++++++++++++++++++++++++++++%

La forme de Maurer-Cartan d'un groupo\"ide de Lie est une $1$-forme sur le groupo\"ide invariant sous l'action du groupo\"ide sur lui-m\^eme
par composition au but. Pour un $\D$-groupo\"ide de Lie, la formes de Maurer-Cartan proviennent de la restriction de celle
de $J^*$. 

\subsubsection{Forme de Maurer-Cartan de $J^*$}

L'espace $J^*_q$ n'agit pas sur son espace tangent et la construction de la forme invariante ne peut se faire qu'en consid\'erant l'action 
de $J^*_{q+1}$ sur le tangent de $J^*_q$. 

Plus pr\'ecisement, la structure de groupo\"ide de $J^*_q$ donn\'ee par la composition
$$
 c : (J^*_q,t) \times_{X} (J^*_q,s) \to J^*_q.
 $$
 se prolonge en une structure de groupo\"ide sur l'espace $J^*_1({J^*_q} _{/source})$ des jets de sections de la projection \textit{source} 
 se projetant sur les sections de $J^*_1(X\times X _{/source})$.
 Ce groupo\"ide agit naturellement sur $T{J^*_q}_{/source}$. 
 Cette action se comprend plus facilement en utilisant la pr\'esentation suivante. Soit $\varphi$ une section analytique locale de ${J^*_q} _{/source}$ sur un ouvert 
 $U$ dont la composante d'ordre $0$, $\varphi_0$ est inversible. Par les formules de composition, elle induit un diff\'eomorphisme de l'ouvert de 
 ${J^*_q}$ des jets de buts dans $U$ sur l'ouvert des jets de but dans $\varphi_0(U)$. Le diff\'eomorphisme induit sur le tangent de 
 ${J^*_q}$ ne d\'epend que du premier jet de la section et donne l'action. 
 
 Nous noterons $Tc$ cette action :  
  $$
      Tc :  (T {J^*_q}_{/source}, t)  \times_{X}  (J^*_1 ({J^*_q}_{/source} ), s) \to T {J^*_q}_{/source}.
  $$
%aaaaaaaa
Cette action induit une action de $J^*_{q+1}$ sur $T {J^*_q}_{/source}$ par l'inclusion canonique $J^*_{q+1} \subset J^*_1({J^*_q}_{/source})$
et donne un isomorphisme 
%aaaaaaaa
      $$
     T{ J^*_q}_{/source}|_{id}  \times_{X}  J^*_{q+1} \to T{J^*_q}_{/source}.
      $$
%aaaaaaaa
Apr\`es 
identification de $T{J^*_q}_{/source}|_{id}$ avec $J_q(TX_{/X})$, l'image de la premi\`ere projection donne une application 
%aaaaaaaa
      $$
      \Theta : T{J^*_q}_{/source} \to J_q(TX_{/X})
      $$
%aaaaaaaa
invariante sous $J^*_{q+1}$. 
%
%----------------%
\begin{definition}
%----------------%
   Cette forme est la forme de Maurer-Cartan de $J^*_{q+1}$.
   
   Les formes construites ainsi \'etant compatibles aux projections $\pi^{q+k}_{q}$, la forme induite sur $J^*$ sera appel\'ee forme de Maurer-Cartan de $J^*$. 

%----------------%
\end{definition}
%----------------%

\begin{remarque}
Fixons un point $x_0 \in X$ et identifions la sous-vari\'et\'e des jets de source $x_0$ avec l'espace des rep\`eres, $R_q(X)$ 
ainsi que $J_q(TX_{/X})|_{x_0}$ avec $ J_q(\chi(\CC^n,0))$ l'espace des jets d'ordre $q$ de champs de vecteurs en $0$. On obtient alors une forme 
sur  $R_q(X)$ \`a valeurs dans $ J_q(\chi(\CC^n,0))$ invariante sous l'action de $J^*_{q+1}$.  Nous noterons
$$
      \Theta|_{x_0} : TR_q(X) \to J_q(\chi(\CC^n,0))
$$
cette forme.
%aaaaaaaa
\end{remarque}

Il existe un second crochet sur $J_q(TX_{/X})$ diff\'erent du crochet de Spencer. Il est  \`a valeurs dans 
$J_{q-1}(TX_{/X})$. C'est le crochet fibre \`a fibre. Il est d\'efini en coordonn\'ees locales par les formules 
%aaaaaaaa
$$
\{u,v\} = j_{q-1}[\widetilde{u}, \widetilde{v}]
$$
%aaaaaaaa
o\`u $u$ et $v$ sont des jets en un 
point $x\in X$, $\widetilde{u}$ et $\widetilde{v}$ sont des champs formels en $x$ les prolongeant, $[\ ,\ ]$ est le crochet de Lie et $j_{q-1}$ est
la prise du jet d'ordre $q-1$.

La forme de Maurer-Cartan v\'erifie des \'equations de structure reliant la diff\'erentielle relative \`a la projection \textit{source}
et le crochet fibre \`a fibre 
sur $J_q(TX_{/X})$ :
%aaaaaaaa
      $$
      d_{/s} \pi_{q-1}^q \Theta = - \frac{1}{2} \{ \Theta , \Theta \},
      $$
%aaaaaaaa
o\`u $\pi_{q-1}^q$ est la projection de $J_q(TX_{/X})$ sur $J_{q-1}(TX_{/X})$.
Nous renvoyons \`a \cite{Guillemin-Sternberg, s-s} pour plus de d\'etails.

%aaaaaaaa
 Ces \'equations ce restreignent aux espaces de jets de source fix\'ee.
En choisissant une base de l'espace vectoriel des jets d'ordre $q$ de champs de vecteurs en $x_0$ et en \'ecrivant $\Theta|_{x_0}$ en
coordonn\'ees, on obtient une famille de 1-formes invariantes sur $R_q(X)$. 
En choisissant la base monomiale $\frac{x^\alpha}{\alpha !} \frac{\partial}{\partial x_i}$, on obtient une suite de formes  $\theta_i^\alpha$ 
satisfaisant les identit\'es diff\'erentielles :
$$
d \theta_i^\alpha = \sum_j \theta_j^0 \wedge \theta_i^{\alpha + \epsilon_j} + \sum_{j,|\beta|\geq 1} \left( \! \begin{array}{c} \alpha \\ \beta \end{array} \! \right) \theta_j^\beta \wedge \theta_i^{\alpha-\beta + \epsilon_j} 
$$
pour $|\alpha| \leq q-1$ et $\left( \! \begin{array}{c} \alpha \\ \beta \end{array} \! \right) = \prod \left( \! \begin{array}{c} \alpha_i \\ \beta_i \end{array} \! \right)$ 
si $\alpha_i \geq \beta_i$ et z\'ero sinon.
Nous les appellerons formes de Cartan d'ordre $q$ et les identit\'es pr\'ec\'edentes \'equations de structures d'ordre $q$.\\

Le groupe alg\'ebrique des jets d'ordre $q+1$ de source et but $x_0$, $\Gamma_{q+1}^n$, agit sur $R_q(X)$ par composition \`a la source. 
La composition au but commutant avec la
composition \`a la source, $\Gamma_{q+1}^n$ transforme des formes invariantes en des formes invariantes satisfaisant les m\^emes identit\'es
diff\'erentielles. Soit $\varphi$ un \'el\'ement de $\Gamma_{q+1}^n$, notons $s_\varphi$ l'action par composition \`a la source sur $R_q(X)$
et $\varphi^*$ celle sur $J_q(\chi(\CC^n,0))$. On a la formule suivante :
$$
s_\varphi^* \Theta|_{x_0} = \varphi^* \circ \Theta|_{x_0}.
$$

\subsubsection{Forme de Maurer-Cartan de $\Y$}

Consid\'erons maintenant $\Y$ un $\D$-groupo\"ide de Lie agissant sur $X$.  
%----------------%
\begin{definition}
%----------------%
   La restriction de $\Theta$ sur $\Y_{q+1}$ est \`a valeurs
   dans le sous-espace vectoriel $L\Y_{q}$ :
   %aaaaaaaa
      $$
      \Theta_{\Y} : T {\Y_q}_{/X} \to L\Y_{q} \subset J_q(TX_{/X}).
      $$
   %aaaaaaaa
   C'est la forme de Maurer-Cartan d'ordre $q$ de $\Y$.
\end{definition}
   
   Dans la suite nous nous restreindrons \`a l'\'etude des $\D$-groupo\"ides de Lie transitif et travaiilerons dans le cadre de la d\'efinition suivante.
   
   \begin{definition}
   Soit $\Y$ un $\D$-groupo\"ide de Lie sur $X$ transitif.
   Fixons une source $x_0$ r\'eguli\`ere
   pour $\Y$ (\textit{i.e.} hors du lieu singulier donn\'ee dans dans la d\'efinition \ref{groupoidedelie}) et notons
   $R_q(\Y)$ l'espace des jets de $\Y$ de source $x_0$. On d\'efini une forme invariante sur $R_q{\Y}$ en restreignant la forme de Maurer-Cartan :
   $$
      \Theta_{\Y}|_{x_0} : TR_q(\Y) \to J_q(\chi(\CC^n,0)).
   $$
%----------------%
\end{definition}
%----------------%   

Les orbites de l'action de $\Y_{q+1}$ sur $R_q(X)$ sont des sous-vari\'et\'es d\'ecrites par le th\'eor\`eme \ref{theoinvdiff}.
Une orbite particuli\`ere est donn\'ee par $R_q(\Y)$ mais le groupe $\Gamma_{q+1}^n$ agissant 
par composition en $x_0$ m\'elangent ces orbites.
La restriction de $\Theta$ sur une orbite permet de construire un syst\`eme de formes invariantes satisfaisant des \'equations de structures particuli\`eres.

\subsubsection{Suites de Godbillon-Vey g\'en\'erales}
%+++++++++++++++++++++++++++++++++++++++++++++++++++++++%
%
\begin{definition}

Une suite de Godbillon-Vey pour un feuilletage de codimension $p$ est une suite de $1$-formes rationnelles 
$\{\omega_i^\alpha ; i \in\{1,\ldots p\} ,\alpha \in \NN^p \}$ telles que $\{\omega_i^0 ; i \in\{1,\ldots p\} \}$ d\'efinissent le feuilletage et
$$
d \omega_i^\alpha = \sum_j \omega_j^0 \wedge \omega_i^{\alpha + \epsilon_j} + \sum_{j,|\beta|\geq 1} \left( \! \begin{array}{c} \alpha \\ \beta \end{array} \! \right) \omega_j^\beta \wedge \omega_i^{\alpha-\beta + \epsilon_j} 
$$
avec $\left( \! \begin{array}{c} \alpha \\ \beta \end{array} \! \right) = \prod \left( \! \begin{array}{c} \alpha_i \\ \beta_i \end{array} \! \right)$ 
si $\alpha_i \geq \beta_i$ et z\'ero sinon.
\end{definition}

\begin{proposition}
Soit $\F$ un feuilletage singulier de $X$.
Quitte \`a se placer sur  $\widetilde{X} \rightarrow X$ une vari\'et\'e finie au-dessus d'un ouvert de Zariski dense de $X$, il existe toujours une
suite de Godbillon-Vey pour $\F$.
\end{proposition}

Cette proposition est classique, nous allons en donner une preuve utilisant les constructions pr\'ec\'edente.
Nous allons d'abord construire
ces formes sur le $\D$-groupo\"ide de Lie $Aut(\F)$ et les descendre ensuite sur $X$.

\begin{preuve}
Un feuilletage $\F$ d\'efinit un $\D$-groupo\"ide de Lie $Aut(\F)$ : le groupo\"ide des transformations locales pr\'eservant le feuilletage.
Pour tout entier
$q$, $\F_q \subset T{Aut(\F)_q}_{/source}|_{id}$ et l'action de $Aut(\F)$ par composition au but permet de prolonger $\F_q$ en un sous-espace 
vectoriel $\widetilde{\F_q}$ de $T{Aut(\F)_q}_{/source}$. 

Soit $x_0$ un point r\'egulier de $\F$ (\textit{i.e.} au voisinage duquel $\F$ est facteur direct). C'est aussi un point r\'egulier de $Aut(\F)$.
 Les jets solutions de $Aut(\F)_q$ de sources $x_0$ forment une sous-vari\'et\'e de $R_q(X)$ not\'e $R_q(\F)$

\begin{definition}
Les formes sur $R_q(\F)$ s'annulant sur $\widetilde{\F_q}$ seront dites transverses \`a $\F$.
\end{definition}

On construit une forme transverse de la mani\`ere suivante. Soit $N_\F$ la $\D$-alg\`ebre de Lie de $Aut(\F)$ et soit 
$\pi : {{N_\F}_q}|_{x_0} \to V_q$ le quotient par ${\F_q}|_{x_0}$. 
%
%---------------%
\begin{definition}
%---------------%
   La projection, 
   $$
   \pi \circ \Theta_{Aut(\F)}|_{x_0} : T R_q(\F) \to V
   $$ est une forme transverse au feuilletage. 
   Elle sera not\'ee $\Theta_\F$
   
%---------------%
\end{definition}
%---------------%
%

Le crochet sur la fibre de $N_\F$ en $x_0$ induit un crochet sur $V$.
La forme transverse v\'erifie les \'equations de structure induites :
$$
d \Theta_\F = -\frac{1}{2} \left\{ \Theta_\F, \Theta_\F \right\}.
$$
Soit $p$ la codimension du feuilletage. L'espace vectoriel $V_q$ s'identifie \`a $J_q(\chi(\CC^p,0))$. 
En \'ecrivant $\Theta_\F|_{x_0}$ dans la base monomiale $\frac{x^\alpha}{\alpha !} \frac{\partial}{\partial x_i}$ sur  $J_q(\chi(\CC^p,0))$, 
on obtient une suite de 1-formes invariantes transverse \`a $\F$ : $\theta_i^{\alpha} ,\ \ 1 \leq i \leq p , \ \ \alpha \in \NN^p,$ satisfaisant les \'equations de structures
d'ordre $q$ 
en dimension $p$.

Soit $\widetilde{X} \rightarrow X$ une vari\'et\'e au-dessus d'un ouvert de Zariski dense de $X$ tel qu'il existe une section
$f : \widetilde{X} \to \widetilde{X} \times_{X} R_q(\F)$ de la projection \textit{but}. Les formes $\theta_i^0$ s'annullent sur $\widetilde{\F}_0$.
Par 
construction ce feuilletage sur $R_0(\F) = X$ est $\F$. La suite $\{\omega_i^\alpha = f*\theta_i^\alpha\}$ est une suite de Godbillon-Vey pour $\F$.  
\end{preuve}

\subsubsection{Suites de Godbillon-Vey sp\'eciales}

%--------------% 
\begin{definition}
%--------------%
%
%---------------%
   Soient $\Y$ un $\D$-groupo\"ide de Lie transitif admissible pour un feuilletage $\F$ et $\widetilde{X} \rightarrow X$ une vari\'et\'e au-dessus d'un ouvert de Zariski dense de $X$ tel qu'il existe une section
$f : \widetilde{X} \to \widetilde{X} \times_{X} R_q(\Y)$ de la projection \textit{but}. Les formes $\omega_i^\alpha = f^*\theta_i^\alpha$ seront 
appell\'e suite de Godbillon-Vey sp\'ecial associ\'ee a $\Y$
%---------------%
\end{definition}
%---------------%
% 
Dans le cas trivial $\Y = Aut(\F)$, nous avons appell\'e la suite de Godbillon-Vey obtenue
est g\'en\'erale. Lorsque $\Y = Gal(\F)$, nous dirons qu'elle est minimale. Les cas interm\'ediaires sont dits sp\'eciaux.  
Dans la suite de cet article, nous regarderons particuli\`erement le cas $X=\CC^3$. On peut toujours choisir une section \`a coefficients alg\'ebriques
 de la projection 
but et construire ainsi une suite de formes \`a coefficients alg\'ebriques sur $\CC^3$.
 
L'objet de cet article est l'\'etude d'un feuilletage de codimension deux.
La classification locale des pseudogroupes de Lie r\'eguliers agissant sur $\CC^2$ a \'et\'e donn\'ee par S. Lie \cite{Lie}. 
Dans \cite{lessousgroupesdesgroupes}, \'E. Cartan donne une preuve de la classification de S. Lie en utilisant ces 
\'equations de structures. Appliqu\'ee \`a la d\'etermination des \'equations de structures possibles pour les formes \`a valeurs dans 
$J_q(\chi(\CC^2,0))$,
la preuve de Cartan permet d'obtenir le r\'esultat suivant.
%
%--------------%
\begin{theoreme}
%--------------%
\label{itchiban}
%--------------%
    Soit $\F$ un feuilletage de $\CC^n$ de codimension deux, d\'efini par une $2$-forme ferm\'ee $\gamma$. Le $\D$-groupo\"ide de Lie $Inv(\gamma)$
    d'invariance de cette forme est un $\D$-groupo\"ide de Lie admissible pour $\F$. Si le groupo\"ide de Galois de $\F$ n'est pas $Inv(\gamma)$ alors
    on est dans un des cas suivant :
    \begin{itemize}
       \item $Gal(\F)$ est intransitif : $\F$ admet une int\'egrale premi\`ere rationnelle,
       \item $Gal(\F)$ est imprimitif en codimension un : il existe une 1-forme \`a coefficients alg\'ebriques int\'egrable s'annulant sur $\F$,
       \item $Gal(\F)$ est transversalement affine : il existe un vecteur de 1-formes \`a coefficients alg\'ebriques   
      \mbox{ $\Omega^0 = \left(  \! \begin{array}{c} \omega^0_1 \\ \omega^0_2 \end{array}  \! \right)$} d\'efinissant le 
       feuilletage et une matrice de 1-formes \`a coefficients alg\'ebriques $\Omega^1$ de trace nulle telle que:
       $$
       d\Omega^0 = \Omega^1 \wedge \Omega^0 \ \ \textrm{ et } \ \  d\Omega^1 = \Omega^1 \wedge \Omega^1.
       $$
       Une telle suite sera appel\'ee suite de Godbillon-Vey de type $asl_2$ ou encore $asl_2$-suite.
    \end{itemize}
%------------% 
\end{theoreme}
%------------%
 Nous donnons en annexe la preuve de Cartan de ce th\'eor\`eme.

\section{Le groupo\"ide de Galois de $P_1$}
 
 Le feuilletage associ\'e \`a la premi\`ere \'equation de Painlev\'e, $\frac{d^2y}{dx^2}=6y^2 + x$, est donn\'ee par le champ de vecteurs 
 $$
 X_1 = \frac{\partial}{\partial x} + y'\frac{\partial}{\partial y} + (6y^2 +x)\frac{\partial}{\partial y'}
 $$
 sur $\CC^3$ avec les coordonn\'ees $(x,y,y')$.
 
 Ce champ de vecteurs \'etant de divergence nulle, le feuilletage est d\'efini par la 2-forme ferm\'ee $\gamma = i_{X_1}dx \wedge dy \wedge dy'$. 
 %
%-------------%
\begin{theoreme}
%-------------%
%
\label{groupoidedegaloisdeP1}
%
%-------------%
   Le groupo\"ide de Galois de $P_1$ est le $\D$-groupo\"ide de Lie  $Inv(\gamma)$ laissant la forme $\gamma$ invariante.
%-------------%
\end{theoreme}
%-------------% 
 
On consid\`erera le champ $X_1$ dans la famille triviale
 $$
 X_\alpha = \frac{\partial}{\partial x} + y'\frac{\partial}{\partial y} + (6y^2 +\alpha^5 x)\frac{\partial}{\partial y'}, \ \ \alpha \not = 0.
 $$
 Ce champ est de poids $-1$ sous l'action de
 $$
 \Sigma = x \frac{\partial}{\partial x} -2y \frac{\partial}{\partial y} -3y' \frac{\partial}{\partial y'}- \alpha \frac{\partial}{\partial \alpha}
 $$
 c'est-\`a-dire que $[\Sigma, X_\alpha]= -X_\alpha$. 
 Le flot de $\Sigma$ donne une \'equivalence orbitale entre les diff\'erents champs $X_\alpha$ tant que $\alpha \not = 0$. Lorsque $\alpha=0$,
 la premi\`ere \'equation de Painlev\'e d\'eg\'en\`ere sur une \'equation particuli\`erement simple.

 L'organisation des calculs se fera suivant les id\'ees de J. Drach \cite{Drach-painleve}. Soit $\pi_\alpha$ la projection alg\'ebrique de 
 $\CC^3 \times \CC_\alpha$ sur $\CC^3$ donn\'ee par le flot de $\Sigma$ :
 $$
    \pi_\alpha (x,y,y',\alpha) = (\alpha x, \frac{y}{\alpha^{2}}, \frac{y'}{\alpha^{3}},1).
 $$ 
 Supposons qu'un des trois syst\`emes d'\'equations aux d\'eriv\'ees partielles donn\'e par le th\'eor\`eme \ref{itchiban} admette une solution alg\'ebrique. En consid\'erant leurs images 
 inverses par $\pi_\alpha$ et en d\'eveloppant les \'equations et la solution en puissances de $\alpha$, on obtient une suite de syst\`emes plus 
 simples ayant encore des solutions alg\'ebriques. Ces \'equations sont de la forme $X_0 R = *$. On les \'etudiera dans les 
 coordonn\'ees adapt\'ees $x,y,u=y'^2-4y^3$.
 
 Nous noterons $p_\Sigma(\ .\ )$ le poids d'une fonction ou d'une forme homog\`ene sous $\Sigma$ 
 et $deg_{a,b,\ldots}$ les degr\'es par rapport
 aux variables $a,b,\ldots$ ayant chacune le degr\'e $1$.   

 Un polyn\^ome $P(x,y,y')$ se prolonge en une fraction  de poids nul sur 
 $\CC^3 \times \CC_\alpha$ par 
 $$
 \pi_\alpha^*P = P(\alpha x, \frac{y}{\alpha^{2}}, \frac{y'}{\alpha^{3}}) \in \CC(x,y,y',\alpha)
 $$
  et en un polyn\^ome en $x$, $y$, $y'$ et $\alpha$ homog\`ene pour $\Sigma$ en multipliant par une puissance convenable de $\alpha$. 
 Nous prolongerons les formes de $\CC^3$, $a\,dx+b\,dy+c\,dy'$, \`a $\CC^3 \times \CC_\alpha$ par $(\pi_\alpha^*a) \alpha dx+(\pi_\alpha^*b) \frac{dy}{\alpha^{2}}+(\pi_\alpha^*c)\frac{dy'}{\alpha^{3}}$. 
 Une forme int\'egrable ne se prolongera pas en une forme int\'egrable mais seulement int\'egrable modulo $d\alpha$.
 
 Le feuilletage donn\'e par le champ $X_1$ est d\'efini par les deux $1$-formes 
 $$
 \omega^0_1 = y'dy' -(6y^2+x)dy,\ \ \ \  \omega^0_2 = \frac{dy}{y'} - dx.
 $$ 
 Nous noterons $\omega^0_1|_\alpha = y'dy' -(6y^2+\alpha^5 x)dy$ et 
 $\omega^0_2|_\alpha = \frac{dy}{y'} - dx$ leurs prolongements au param\`etre $\alpha$ et $\omega^0_1|_0$, $\omega^0_2|_0$ les restrictions 
 \`a $\{\alpha = 0 \}$. Ces formes forment les vecteurs $\Omega^0$, $\Omega^0|_\alpha$ et $\Omega^0|_0$.

%%%%%%%%%%%%%%%%%%%%%%%%%%%%%%%%%%%%%%%%%%%%%%%%%%%%%%%%%%%%%%%%%%%%%%%%%%%%%%%%%%%%%%%%%%%%%%%%%%%%%%%%%%%%%%%%%%%%%%%%%%%%
%%%%%%
%%                 ETUDE DE X_0 
%%                 ETUDE DE X_0 
%%                 ETUDE DE X_0 
%%                 ETUDE DE X_0 
%%                 ETUDE DE X_0 
%%                 ETUDE DE X_0 
%%                 ETUDE DE X_0 
%%                 ETUDE DE X_0 
%%                 ETUDE DE X_0 
%%                 ETUDE DE X_0 
%%                 ETUDE DE X_0 
%%                 ETUDE DE X_0 
%%                 ETUDE DE X_0 
%%                 ETUDE DE X_0 
%%                 ETUDE DE X_0 
%%                 ETUDE DE X_0 
%%                 ETUDE DE X_0 
%%                 ETUDE DE X_0 
%%                 ETUDE DE X_0 
%%                 ETUDE DE X_0 
%%%%%%
%%%%%%%%%%%%%%%%%%%%%%%%%%%%%%%%%%%%%%%%%%%%%%%%%%%%%%%%%%%%%%%%%%%%%%%%%%%%%%%%%%%%%%%%%%%%%%%%%%%%%%%%%%%%%%%%%%%%%%%%%%%%

 \subsection{\'Etude pr\'eliminaire de l'\'equation r\'eduite $X_0$}
%-------------------------------------------------------------------%

 Le champ de vecteurs 
 $$
 X_0 = \frac{\partial}{\partial x} + y'\frac{\partial}{\partial y} + 6y^2\frac{\partial}{\partial y'}
 $$
 admet une int\'egrale premi\`ere rationnelle $u=y'^2 -4y^3$ et une autre transcendante $x-\int \frac{1}{4y^3+u}dy$.
  Nous utiliserons souvent les formules suivantes :
    %aaaaaaaaa
    \begin{equation}
    \label{formules}
    \begin{array}{rcl}
      X_0 \left( \left( x+2\frac{y}{y'} \right)^{n+1} \right) &=& \left( x+2\frac{y}{y'} \right)^n \frac{(n+1)3u}{y'^2} \\
      X_0 \left( \left( 7x+20\frac{y}{y'}-24\frac{y^4}{y'^3} \right)^{n+1} \right) &=& \left( 7x+20\frac{y}{y'}-24\frac{y^4}{y'^3} \right)^n \frac{(n+1)27u^2}{y'^4} . 
    \end{array}
    \end{equation}
    %aaaaaaaaa

\begin{lemme}
%-------------%
  \label{prim0}
%-------------%
  Si $\eta_0$ est une 1-forme polynomiale, int\'egrable, telle que $\eta_0(X_0) = 0$ et $i_{X_0}d\eta_0 = 0$ alors $\eta_0 = f(u)du$ ou 
  $\eta_0 = f(u) \left(udx -(2xy^2 + y')dy + \frac{1}{3}(xy'+2y)dy' \right)$
%-------------%
\end{lemme}

\begin{preuve}
   \'Ecrivons $\eta_0 = a\,du + b\,(dy-y'dx)$. Si $b=0$, on trouve la premi\`ere expression. Sinon, en \'ecrivant la condition d'int\'egrabilit\'e de $\eta_0$, on 
   obtient l'\'equation
   $$
   -X_0\left( \frac{a}{b}\right) + \frac{6y^2}{y'} \frac{a}{b} - \frac{1}{2y'} =0.
   $$ 
   En utilisant les formules (\ref{formules}), on int\`egre cette \'equation. On trouve 
   $$
   \frac{a}{b} = -\frac{1}{6u} (xy' +2y).
   $$
   On en d\'eduit que $\eta_0 = f \left( (\frac{x+2y/y'}{6u}) du - \frac{1}{y'} (dy-y'dx)\right )$. En \'ecrivant la condition $i_{X_0}d\eta_0 =0$, 
   on trouve que $X_0f =0$, donc $f$ est une fonction de $u$.
   En d\'eveloppant $\eta_0$, on trouve la seconde expression. 
%-------------%  
\end{preuve}
%-------------%
%
%
%-------------%
\begin{lemme}
%-------------%
   \label{asl20}
%-------------% 
   Le feuilletage d\'efini par $X_0$ n'admet pas de suite de Godbillon-Vey rationnelle de type $asl_2$ commen\c{c}ant par $\omega^0_1|_0$ et 
   et $\omega^0_2|_0$.
%-------------%
\end{lemme}
%-------------%
%
%
%-------------%
\begin{preuve}
%-------------%
%
   Supposons qu'il existe une matrice de 1-formes $\Omega^1|_0$ de trace nulle satisfaisant :
   %aaaaaaaaa
   $$
    \begin{array}{rcl}
      d\Omega^0|_0 &=& \Omega^1|_0 \wedge \Omega^0|_0 \\
      d\Omega^1|_0 &=& \Omega^1|_0 \wedge \Omega^1|_0 .
    \end{array}
   $$
   %aaaaaaaaa
   On peut supposer en toute g\'en\'eralit\'e que ces formes sont homog\`enes sous 
   $\Sigma$.
   La premi\`ere \'egalit\'e donne 
   %aaaaaaaaa
   $$
    \Omega^1|_0 = \left( \begin{array}{cc} 
                     0             & 0 \\
                     \frac{1}{y'^3}& 0  
                    \end{array}           \right)dy + A\, \omega^0_1|_0 + B\, \omega^0_2|_0
   $$
   %aaaaaaaaa
   avec $A\left( \begin{array}{c}
                          0\\
                          1 \end{array} \right) =B\left( \begin{array}{c}
                                                                 1\\
                                                                 0 \end{array} \right) $ et $trace\, A = trace\, B = 0$.\\
   La deuxi\`eme \'egalit\'e donne un premier syst\`eme d'\'equations que nous allons r\'esoudre : 
   %aaaaaaaa
   $$
    i_{X_0} d\Omega^1|_0 = \left[\Omega^1|_0(X_0), \Omega^1|_0 \right].
   $$
   %aaaaaaaa
   Calculons les deux membres de cette \'egalit\'e :\\
    
    $\bullet \ i_{X_0} d\Omega^1|_0 = \left( \left( \begin{array}{cc} 0 & 0 \\ 1/y'^4 & 0  \end{array} \right) + X_0 A + \frac{1}{y'^2} B \right) \omega^0_1|_0 
                                                 +  X_0B\, \omega^0_2|_0$ \\
                                                 
    $\bullet  \left[ \Omega^1|_0(X_0),\Omega^1|_0\right] = \left[ \left( \begin{array}{cc} 0 & 0 \\ 1/y'^2 & 0  \end{array} \right), A \right] \omega^0_1|_0 
                                                                              + \left[ \left( \begin{array}{cc} 0 & 0 \\ 1/y'^2 & 0  \end{array} \right),B\right] \omega^0_2|_0$ \\
                                                                              
    \noindent                                  
    En identifiant les coefficients de $\omega^0_1|_0$ et ceux de $\omega^0_2|_0$, on obtient les deux syst\`emes suivants,
    %aaaaaaaaa
    \begin{subeqnarray}
    X_0 A &=& \left[ \left( \begin{array}{cc} 
                     0             & 0 \\
                     \frac{1}{y'^2}& 0  
                    \end{array}           \right),  A \right] - \left( \begin{array}{cc} 
                                                                                0             & 0 \\
                                                                                \frac{3}{y'^4}& 0  
                                                                        \end{array}           \right) - \frac{1}{y'^2}B  \slabel{asl20-1}\\
    X_0 B &=& \left[ \left( \begin{array}{cc} 
                     0             & 0 \\
                     \frac{1}{y'^2}& 0  
                    \end{array}           \right),  B \right]. \slabel{asl20-2}
    \end{subeqnarray}                                                                             
    %aaaaaaaaa 
    Si $B=\left( \begin{array}{cc} 
                     a & b \\
                     c & -a  
                 \end{array}   \right)$,
    on \'ecrit puis r\'esoud le syst\`eme (\ref{asl20-2}) en utilisant les formules (\ref{formules}) : 
    %aaaaaaaaaa                                                                                       
    \begin{equation}
    \label{labas}
    \begin{array}{rclcrcl}
      X_0 b &=&0               & \; \; \Rightarrow \; \; & b&=&b(u), \\
      X_0 a &=& \frac{-b}{y'^2}& \; \; \Rightarrow \; \; & a&=&\frac{-b}{3u}( x + 2y/y'), \\
      X_0 c &=& \frac{2a}{y'^2}& \; \; \Rightarrow \; \; & c&=&\frac{-b}{9u^2}( x + 2y/y')^2,
    \end{array}
    \end{equation}
    %aaaaaaaaaa
    sous l'hypoth\`ese $b \not = 0$ (les r\'esultats sont analogues lorsque $b$ est nul).
    On r\'esoud ensuite l'\'equation (\ref{asl20-1}). Si $A=\left( \begin{array}{cc} 
                                                                     -c & a \\
                                                                     d & c  
                                                           \end{array}   \right)$,
    le syst\`eme donne 
    %aaaaaaaaaaaaa%
    $$
      X_0 d =\frac{-3c}{y'^2} - 3\frac{1}{y'^4}  \  \Rightarrow \  d= \frac{b}{27u^3}( x + 2y/y')^3 - \frac{1}{9u^2} \left( 7x+20\frac{y}{y'}-24\frac{y^4}{y'^3} \right).\\
    $$
    %aaaaaaaaaaaaa%
    On obtient une contradiction en consid\'erant l'\'equation suivante :
    %aaaaaaaaaaaaa%
    $$
    i_{\frac{\partial}{\partial x}}d\Omega^1|_0 = \left[ \Omega^1|_0\left(\frac{\partial}{\partial x}\right) , \Omega^1|_0 \right].
    $$
    %aaaaaaaaaaaaa%
    En identifiant les coefficients de $\omega^0_1|_0$ dans cette \'equation, on obtient
    %aaaaaaaaaaaaa%
    $$
    \frac{1}{y'}\frac{\partial B}{\partial y'} + \frac{\partial A}{\partial x} = \left[ A,B \right],
    $$
    %aaaaaaaaaaaaa%
    c'est-\`a-dire,
    %aaaaaaaaaaaaa%
    \begin{subeqnarray}
    \frac{1}{y'}\frac{\partial a}{\partial y'} + \frac{\partial c}{\partial x} &=& ac-db \\
    \frac{1}{y'}\frac{\partial b}{\partial y'} + \frac{\partial a}{\partial x} &=& - 2(bc+a^2), \slabel{dessus}\\
    \frac{1}{y'}\frac{\partial c}{\partial y'} + \frac{\partial d}{\partial x} &=& 2(da-c^2). \slabel{dessus'}
    \end{subeqnarray}
    %aaaaaaaaaaaaa%
    L'\'equation (\ref{dessus})  s'\'ecrit $2b'(u)-\frac{1}{3u}b(u)=0$. La seule solution rationnelle est $b=0$. 
    Dans ce cas $a=a(u)$ et l'\'equation (\ref{dessus}) donne $a=0$. Les \'equations (\ref{labas}) impliquent $c=c(u)$ et 
    $d=\frac{-3c}{y'^2}-\frac{1}{9u^3}\left( 7x+20\frac{y}{y'}-24\frac{y^4}{y'^3} \right)$. L'\'equation (\ref{dessus'})
    se r\'e\'ecrit alors
    %aaaaaaaaaaaaa% 
    $$
    2c'(u)-\frac{7}{9u^3}=-2c^2.
    $$
    %aaaaaaaaaaaaa%
    Une solution de cette \'equation ne pouvant pas avoir de p\^oles d'ordre entier en 0, elle ne peut pas avoir de solution rationnelle. 
    On aboutit \`a une contradiction qui prouve le lemme.
\end{preuve}

%%%%%%%%%%%%%%%%%%%%%%%%%%%%%%%%%%%%%%%%%%%%%%%%%%%%%%%%%%%%%%%%%%%%%%%%%%%%%%%%%%%%%%%%%%%%%%%%%%%%%%%%%%%%%%%%%%%%%%%%%%%%
%%%%%%%%
%%        TRANSITIVITE
%%        TRANSITIVITE
%%        TRANSITIVITE
%%        TRANSITIVITE
%%        TRANSITIVITE
%%        TRANSITIVITE
%%        TRANSITIVITE
%%        TRANSITIVITE
%%        TRANSITIVITE
%%        TRANSITIVITE
%%        TRANSITIVITE
%%        TRANSITIVITE
%%        TRANSITIVITE
%%        TRANSITIVITE
%%        TRANSITIVITE
%%        TRANSITIVITE
%%        TRANSITIVITE
%%        TRANSITIVITE
%%        TRANSITIVITE
%%        TRANSITIVITE
%%        TRANSITIVITE
%%        TRANSITIVITE
%%        TRANSITIVITE
%%%%%%%%
%%%%%%%%%%%%%%%%%%%%%%%%%%%%%%%%%%%%%%%%%%%%%%%%%%%%%%%%%%%%%%%%%%%%%%%%%%%%%%%%%%%%%%%%%%%%%%%%%%%%%%%%%%%%%%%%%%%%%%%%%%%
 
 \subsection{Le groupo\"ide de Galois de $X_1$ est transitif}
%------------------------------------------------------------%
 
Le r\'esultat suivant provient de \cite{Painleve-bulletin}, on pourra aussi consulter \cite{Drach-painleve, Ince}. 
 
 \begin{proposition}
%*************%   
   \label{transitif}
%*************% 
   Le champ $X_1$ n'a pas d'int\'egrale premi\`ere rationnelle.
%*************%
 \end{proposition}

 \begin{preuve}
     Si $H_1$ est une telle int\'egrale, prolongeons-la en une int\'egrale alg\'ebrique du champ $X_\alpha$. On peut \'ecrire $H_\alpha = H_0 + H_1 \alpha + \ldots$ 
     apr\`es avoir fait une division suivant les puissances 
     croissantes et multipli\'e par une puissance convenable de $\alpha$. On peut supposer que $H_0$ n'est pas constante.
     En d\'eveloppant l'\'equation $X_\alpha H_\alpha = 0$ suivant les puissances de $\alpha$, on obtient pour $0\leq i \leq4$ 
     $$
     X_i H_i = 0,
     $$ 
     qui implique $H_i = f_i(u)$ o\`u $u=y'^2-4y^3$ et $f_i$ est rationnelle. Le corfficient de $\alpha^5$ donne 
     $$
     X_0 H_5 + x \frac{\partial f_0(u)}{\partial y'} = 0.
     $$
     Dans les coordonn\'ees $x,y,u$, cette derni\`ere \'equation s'\'ecrit 
     $$
     \left( \frac{\partial}{\partial x} + \sqrt{4y^3+u}\frac{\partial}{\partial y}\right) \left( H^0_5 + \sqrt{4y^3+u} \, H^1_5 \right) = -2xf'_0(u) \sqrt{4y^3+u}
     $$
     o\`u $H_1 = H^0_5 + \sqrt{4y^3+u}\,H^1_5 $ et $H^i_5 \in \mathbb{C}(x,y,u)$. En identifiant les termes libres de radicaux des deux membres, d'une part, 
     et ceux contenant le radical, d'autre part, on obtient les deux \'equations suivantes : 
 
     %aaaaaaaaaaaaaaa%
     \begin{subeqnarray}
     \frac{\partial H_5^0}{\partial y} + \frac{\partial H_5^1}{\partial x}  &=& -2xf'_0(u) \slabel{X1}  \\
     \frac{\partial H_5^0}{\partial x} + (4y^3 + u)\frac{\partial H_5^1}{\partial y} +6y^2 H_5^1  &=& 0. \slabel{X2}   
     \end{subeqnarray}
     %aaaaaaaaaaaaaaa% 
     En d\'erivant (\ref{X2}) par rapport \`a $x$ et en utilisant (\ref{X1}), on obtient 
     %aaaaaaaaaaaaaaa%
     \begin{equation}
       \label{DX2dx}
     \frac{\partial^2 H_5^0}{\partial x^2} - (4y^3 + u)\frac{\partial^2 H_5^0}{\partial y^2} - 6y^2 \frac{\partial H_5^0}{\partial y}= 0.  \\
     \end{equation}
     %aaaaaaaaaaaaaa%
     Supposons que $H_5^0$ ait un p\^ole d'ordre $n$ par rapport \`a $y$ en $a \in \overline{\CC[x,u]}^{alg}$, la cl\^oture alg\'ebrique de 
     $\CC[x,u]$. 
     En d\'eveloppant $H_5^0$ en s\'erie autour de $a$, on obtient une s\'erie d'\'equations diff\'erentielles. La premi\`ere,
     %aaaaaaaaaaaaaa
     $$
     \left( \frac{\partial a}{\partial x} \right)^2 - (4a^3+u)=0,
     $$
     %aaaaaaaaaaaaaa
     implique $4a^3 = -u$. La suivante donne $12(n+1)=-6a^2$. $H_5^0$ est donc un polyn\^ome.
     En \'ecrivant $H_5^0 = \alpha_n y^n + \ldots$ dans (\ref{DX2dx}), on obtient $n =0$ puis
     %aaaaaaaaaaaaaa
     $$
     H_1^0 = a(u)x + b(u).
     $$
     %aaaaaaaaaaaaaa
     L'\'equation (\ref{X1}) donne $\frac{\partial H^1_5}{\partial x} = -2xf'_0(u)$. En d\'erivant (\ref{X2}) par rapport \`a $x$, on 
     obtient $\frac{\partial H^1_5}{\partial x} = 0$, ce qui implique $f'_0=0$. Ceci contredit le choix de $H$ et prouve la proposition.
\end{preuve}

\begin{corollaire}
%^^^^^^^^^^^^^^^%
  \label{hypeinv}
%^^^^^^^^^^^^^^^%   
  Le champ $X_1$ n'admet pas d'hypersurface alg\'ebrique invariante.
%^^^^^^^^^^^^^^^% 
   \end{corollaire}
  
\begin{preuve}
     Soit $P$ l'\'equation d'une hypersurface invariante. Prolongeons $P$ au param\`etre $\alpha$. On obtient un polyn\^ome satisfaisant 
     %aaaaaaaaaaa
        $$
        X_\alpha P = L P
        $$
     %aaaaaaaaaaa
     On a $p_\Sigma(L) = -1$ d'une part et d'autre part comme $deg_{\{x,y,y'\}}(X_\alpha) = 1$, $deg_{\{x,y,y'\}}(L) \leq 1$ 
     et un calcul direct donne $deg_{\{x\}}(L)\leq 0$.  Tout ceci implique que $L$ est nul donc que $P$ est constant.
\end{preuve}

%%%%%%%%%%%%%%%%%%%%%%%%%%%%%%%%%%%%%%%%%%%%%%%%%%%%%%%%%%%%%%%%%%%%%%%%%%%%%%%%%%%%%%%%%%%%%%%%%%%%%%%%%%%%%%%%%%%%%%%%%%%%
%%%%%%%%
%%         PRIMITIF
%%         PRIMITIF
%%         PRIMITIF
%%         PRIMITIF
%%         PRIMITIF
%%         PRIMITIF
%%         PRIMITIF
%%         PRIMITIF
%%         PRIMITIF
%%         PRIMITIF
%%         PRIMITIF
%%         PRIMITIF
%%         PRIMITIF
%%         PRIMITIF
%%         PRIMITIF
%%         PRIMITIF
%%         PRIMITIF
%%         PRIMITIF
%%         PRIMITIF
%%         PRIMITIF
%%         PRIMITIF
%%         PRIMITIF
%%         PRIMITIF
%%         PRIMITIF
%%%%%%%%
%%%%%%%%%%%%%%%%%%%%%%%%%%%%%%%%%%%%%%%%%%%%%%%%%%%%%%%%%%%%%%%%%%%%%%%%%%%%%%%%%%%%%%%%%%%%%%%%%%%%%%%%%%%%%%%%%%%%%%%%%%%%

\subsection{Le groupo\"ide de Galois de $X_1$ est transversalement primitif}

% \begin{preuve}
%      On peut supposer $\omega{\frac{\partial}{\partial x}} = 1$.
%      Soient $\Theta_1 = y'dy' - (6y^2 + x)dy$ et $\Theta_2 = \frac{dy}{y'}-dx$ deux formes d\'efinissant le feuilletage et 
%      $\omega_1,\ldots, \omega_n$ les 1-formes conjugu\'ees de $\omega$. Le champ $X_1$ \'etant polynomial, $\omega_i (X_1) = 0$ 
%      et $\omega_i \wedge \omega_j = a_{ij} \Theta_1 \wedge \Theta_2$. 
%      Consid\'erons la forme $w$ donn\'ee par $\frac{\sum_i \omega_i }{\sqrt{A}}$ si $A = \sum_{i>j} a_{ij}$ est non nul et par
%      $\sum_i \omega_i$ sinon. Prouvons que cette forme est int\'egrable :
%      %aaaaaaaaaaa
%      $$
%      \begin{array}{rcl}
%      d(\frac{\sum_i \omega_i }{\sqrt{A}})\wedge \frac{\sum_j \omega_j }{\sqrt{A}} & = & \sum_{i,j} d\left( \frac{\omega_i}{\sqrt{A}}\right)\wedge\left( \frac{\omega_i}{\sqrt{A}}\right)\\
%                                                                           &=&\sum_{i>j} d \left(\frac{\omega_i\wedge \omega_j}{A}\right)\\
%                                                                           &=&d(\Theta_1\wedge\Theta_2)\\
%                                                                           &=&0.
%      \end{array}
%      $$
%      %aaaaaaaaaaa
%      Le m\^eme calcul montre l'int\'egrabilit\'e dans le cas $A$ nul.
%      Le forme $\sum_i \omega_i$ est une 1-forme rationelle non nulle ($\sum_i \omega_i (\frac{\partial}{\partial x}) = n$) 
%      satisfaisant les hypoth\`eses du lemme. En la multipliant par un facteur convenable on prouve le lemme.  
%      \end{preuve} 
%  
%      

\begin{proposition}
%-------------%   
   \label{prim}
%-------------%
   Le feuilletage donn\'e par $X_1$ n'est inclus dans aucun feuilletage de codimension un donn\'e par une 1-forme alg\'ebrique.
%-------------%
\end{proposition}

     Supposons qu'il existe une 1-forme alg\'ebrique $\eta$ int\'egrable, $\eta \wedge d\eta =0$ et telle que 
     $\eta(X_1)=0$. 
     
\begin{lemme}
%-------------%
  \label{rationnalitefeuilletage}
%-------------%
   Une telle forme peut \^etre choisie polynomiale. 
%-------------% 
\end{lemme}
 
\begin{preuve}
   Soit $\eta$ une telle forme que nous normalisons $\eta = dx+a\,dy+b\,dy'$. Soit $Z$ le lieu de ramification de cette forme sur $\CC^3$. 
   D'apr\`es le corollaire \ref{hypeinv}, $Z$ est transverse aux trajectoires de $X_1$. Pla\c{c}ons-nous au voisinage analytique d'un point de $Z$. 
   N'importe quelle int\'egrale premi\`ere au voisinage d'un point de la trajectoire passant par $p$ se prolonge de mani\`ere univalu\'ee sur un 
   voisinage de $p$. La forme $\eta$ est donc univalu\'ee sur ce voisinage. Le lieu de ramification $Z$ est vide et la forme est rationnelle.
\end{preuve}
         
     Choisissons une de ces formes et prolongeons-la au param\`etre $\alpha$.

\begin{lemme}[\cite{Drach-painleve}]
%------------%
   La forme $\eta_\alpha$  peut \^etre choisie telle que \mbox{$i_{X_\alpha}d \eta_\alpha = 0 $}.
%------------%
\end{lemme}

\begin{preuve}
     Choisissons
     %aaaaaaaaaa
        $$
        \eta_\alpha = -(Py' + Q(6y^2+\alpha^5 x))dx + Pdy + Q dy'
        $$
     %aaaaaaaaaa
     avec $P$ et $Q$ premiers entre eux. 
     En \'ecrivant la condition d'int\'egrabilit\'e, on obtient :
     %aaaaaaaaaa 
        $$
        Q XP - PXQ + 12 yQ^2-P^2 = 0.
        $$
     %aaaaaaaaaa 
     D'apr\`es le th\'eor\`eme de Bezout, il existe donc un polyn\^ome $L$ tel que :
     %aaaaaaaaaa
        $$
        \left\{
        \begin{array}{rl}
           XP+12yQ &=  LP \\
           XQ+P   &=  LQ
        \end{array}
        \right.
        $$
     %aaaaaaaaaa

     En calculant les degr\'es, on obtient que $L$ doit \^etre de degr\'e 1 en $x,y,y'$ et de poids -1 par rapport \`a $\Sigma$.  
     Donc $L= c \alpha ^{2} x$ o\`u $c$ est une constante. En calculant les degr\'es en $x$, on obtient $c=0$. 
     Un calcul direct montre que ceci implique le lemme.
\end{preuve}

\begin{preuved}{de la proposition \ref{prim}}
     Quitte \`a multiplier par un facteur convenable, 
     nous pouvons supposer que $\eta_\alpha$ s'\'ecrit $\eta_0 + \alpha \eta_1 $ + \ldots avec $\eta_i$ des formes polynomiales et $\eta_0 \not = 0$.
     Les \'equations 
     \begin{equation}
     \label{equationeta}
     \begin{array}{rcl}
      \eta_\alpha(X_\alpha) &=& 0 \\
      i_{X_\alpha}d\eta_\alpha & =& 0 
     \end{array}
     \end{equation}
     relient les formes $\eta_0$ et $\eta_1$. \'Ecrivons la forme $\eta_1$ de la mani\`ere suivante
     $$
     \eta_1 = a\, dx + b\,\omega^0_1|_0 +c\,\omega^1_2|_0.
     $$
     Le terme d'ordre $0$ dans les \'equations (\ref{equationeta}) donne $\eta_0(X_0) = 0$ et $i_{X_0}d\eta_0 = 0$. Ces \'equations sont r\'esolues
     par le lemme \ref{prim0}. Celles d'ordre $1$ donne 
     $\eta_0(x\frac{\partial}{\partial y'}) + \eta_1(X_0) = 0$ et $i_{x\frac{\partial}{\partial y'}}d\eta_0 + i_{X_0}d\eta_1 = 0$. 
     En \'ecrivant toutes ces \'equations dans la base $(dx,\omega^0_1|_0,\omega^1_2|_0)$, on obtient  
     \begin{equation}
     \label{equa}
     \begin{array}{rcl}
     a & = & -x \, \eta_0(\frac{\partial}{\partial y'}) \\
     i_{x\frac{\partial}{\partial y'}}d \eta_0 &=& g \, \omega^0_2|_0 \\
     X_0c & = & X_0 a - \frac{\partial a}{\partial x} - g  \\
     X_0b & = & - \frac{c}{y'^2} + \frac{1}{y'} \frac{\partial a}{\partial y'} 
     \end{array} 
     \end{equation}
     Dans chaques cas du lemme \ref{prim0}, nous allons prouver qu'il n'existe pas de triplet $(a,b,c)$ solution de ces \'equations tel que la
     forme $\eta_1$ soit polynomiale.\\
     
     \noindent
     \underline{Premier cas, $p_\Sigma \eta_\alpha \equiv 0 \mod 6$.}\\
     
     D'apr\`es le Lemme \ref{prim0}, $\eta_0 = f(u)du$. Les \'equations portant sur les coefficients de $\eta_1$ deviennent  
     %aaaaaaaaaa
        \begin{subeqnarray}
        a & = & - 2f(u)y'x \slabel {equa-a'}\\
        g & = & 0 \\
        X_0c & = & -12f(u)xy^2 \slabel{equa-c'} \\
        X_0b & = &  -\frac{c}{y'^2} - \frac{2f(u)x}{y'} - 4f'(u)xy' \slabel{equa-b'} 
        \end{subeqnarray}
     %aaaaaaaaaa
     On r\'esoud (\ref{equa-c'}) : $c = 2f(u)(y-xy')$. On obtient ensuite pour l'\'equation (\ref{equa-b'}) :
     %aaaaaaaaaa
        $$
        X_0 b = -4f'(u)xy' - \frac{2f(u)y}{4y^3+u}.
        $$
     %aaaaaaaaaa
     Dans les coordonn\'ees $x,y,u$, o\`u $b=b_0 + \frac{1}{\sqrt{4y^3+u}}b_1$ avec $b_0$ et $b_1$ des polyn\^omes et
     $X_0=  \frac{\partial }{ \partial x}+ \sqrt{4y^3+u}\frac{\partial }{ \partial y}$, l'\'equation ci-dessus se r\'e\'ecrit:
     %aaaaaaaaaa
     \begin{subeqnarray}
        \frac{\partial b_0}{ \partial x} + \frac{\partial b_1}{ \partial y} - \frac{6y^2}{4y^3+u}b_1 &=& - \frac{2f(u)y}{4y^3+u} \slabel{P1}\\
        \frac{\partial b_1}{\partial x} + (4y^3+u)\frac{\partial b_0}{\partial y} &=& - 4f'(u)x (4y^3+u). \slabel{P2}
     \end{subeqnarray}
     %aaaaaaaaaa
     En d\'erivant (\ref{P1}) par rapport \`a $x$ et en tenant compte de la d\'eriv\'ee de (\ref{P2}) par rapport \`a $y$, 
     on obtient
     %aaaaaaaaaa
        \begin{equation}
           \frac{\partial^2 b_0}{\partial x^2 } - (4y^3 + u)\frac{\partial^2 b_0}{\partial y^2 }-6y^2\frac{\partial b_0}{\partial y }=24 f'(u)y^2x.
        \end{equation}
     %aaaaaaaaaa
     On en d\'eduit que $b_0$ est un polyn\^ome d'ordre 1 en $y$ puis 
     %aaaaaaaaaa
        $$
        b_0=4f'(u)xy +a(u)x + b(u) \  ; \   \frac{\partial b_1}{\partial x}=0.
        $$
     %aaaaaaaaaa
     En reportant $b_0$ dans (\ref{P1}), on obtient
     %aaaaaaaaaa
        $$
        (4y^3 + u)\frac{\partial b_1}{\partial y} -6y^2b_1= -16f'(u)y^4-4a(u)y^3 + (f(u)-4uf'(u))y - ua(u),
        $$
     %aaaaaaaaaa
     ce qui implique que le degr\'e de $b_1$ en $y$ est inf\'erieur \`a deux. 
     En posant $b_1=\alpha_2 y^2 + \alpha_1 y^1 +\alpha_0$, on obtient le syst\`eme suivant
     %aaaaaaaaaa
        $$
        \left\{
          \begin{array}{rcl}
             -16f'(u) &=&2\alpha_2\\
               -4a(u) &=&-2\alpha_1\\
                      0 &=& -6\alpha_0\\
         f(u)-4uf'(u) &=&2\alpha_2 u\\
                 ua(u) &=&u\alpha_1.
          \end{array}
        \right.
        $$
     %aaaaaaaaaa
     Ce syst\`eme implique l'existence d'une constante $c$ telle que $f=\frac{c}{u^{1/12}}$. Or $f$ est un polyn\^ome 
     donc $f$ est nul, ce qui contredit le choix de la forme $\eta_\alpha$.\\
     
     \noindent
     \underline{Deuxi\`eme cas, $p_{\Sigma} \eta_\alpha \equiv 1 \mod 6$. }\\
     
     D'apr\`es le lemme \ref{prim0}, $\eta_0 = f(u) \left( \frac{1}{3}\left( x+2\frac{y}{y'}\right)\omega^0_1|_0 - u\, \omega^0_2|_0 \right)$.
     Les \'equations (\ref{equa}) donnent
     %aaaaaaaaaa
        \begin{subeqnarray}
        a & = & - \frac{f(u)}{3}(y'x^2 + 2yx) \slabel {equa-a''}\\
        g & = & -2f(u)uy' - \frac{5f(u)}{3}y' \\
        X_0c & = & -2f(u)y^2x^2  - \left( \frac{7}{3}f(u) - 2f'(u)u\right)y' x \slabel{equa-c''} \\
        X_0b & = &  -\frac{c}{y'^2} -  \frac{f(u)}{3}\frac{x^2}{y'} - \frac{2f'(u)}{3}y'x^2 \slabel{equa-b''} 
        \end{subeqnarray}
     %aaaaaaaaaa
     \'Etudions l'\'equation (\ref{equa-c''}). Le second membre \'etant d'ordre $2$ en $x$, $c$ est d'ordre au plus $3$. 
     Par homog\'en\'eit\'e, $c$ est en fait d'ordre $2$.
     Posons $c = c_2 x^2 +c_1 x+ c_0$. On obtient le syst\`eme suivant 
     $$
     \left\{
     \begin{array}{rcl}
     X_0 c_2 & = & -2f(u) y^2 \\
     X_0 c_1 & = & - \left( \frac{7}{3}f(u) - 2f'(u)u\right)y' -2c_2 \\
     X_0 c_0 & = & -c_1.
     \end{array}
     \right.
     $$
      On les r\'esoud successivement, $c_2 = \frac{-1}{3} f(u) y'$, $c_1 = - \left( \frac{5}{3}f(u) - 2f'(u)u\right)y$ et $X_0 c_0 =\left( \frac{5}{3}f(u) - 2f'(u)u\right)y$.
      Cette derni\`ere \'equation n'a pas de solution polynomiale comme le montre le lemme suivant. Ceci prouve que $f$ doit \^etre nulle, en contradiction avec le choix de $\eta_\alpha$.
      Ceci prouve la proposition.
\end{preuved}
\begin{lemme}
Il n'existe pas de polyn\^ome $R$ satisfaisant $X_0 R = y$.
\end{lemme}
\begin{preuve}
Supposons qu'il existe un tel polyn\^ome. 
Comme le second membre est d'ordre $0$ en $x$, $R$ s'\'ecrit $\alpha_1(u) x + \alpha_0(y,y')$.  En \'ecrivant l'\'equation dans les coordonn\'ees $(x,y,u)$ on a 
$$
\begin{array}{rcl}
     \frac{\partial R^0}{\partial y} + \frac{\partial R^1}{\partial x}  &=& 0  \\
     \frac{\partial R^0}{\partial x} + (4y^3 + u)\frac{\partial R^1}{\partial y} +6y^2 R^1  &=& y.   
\end{array}
$$
Ces \'equations impliquent $(4y^3 + u)\frac{\partial R^1}{\partial y} +6y^2 R^1  = y-\alpha_1(u)$. Cette \'equation n'a pas de solution polynomiale. 
\end{preuve}

%%%%%%%%%%%%%%%%%%%%%%%%%%%%%%%%%%%%%%%%%%%%%%%%%%%%%%%%%%%%%%%%%%%%%%%%%%%%%%%%%%%%%%%%%%%%%%%%%%%%%%%%%%%%%%%%%%%%%%%%%%%%
%%%%%%
%%       ASL2
%%       ASL2
%%       ASL2
%%       ASL2
%%       ASL2
%%       ASL2
%%       ASL2
%%       ASL2
%%       ASL2
%%       ASL2
%%       ASL2
%%       ASL2
%%       ASL2
%%       ASL2
%%       ASL2
%%       ASL2
%%       ASL2
%%       ASL2
%%       ASL2
%%       ASL2
%%       ASL2
%%       ASL2
%%       ASL2
%%       ASL2
%%%%%%
%%%%%%%%%%%%%%%%%%%%%%%%%%%%%%%%%%%%%%%%%%%%%%%%%%%%%%%%%%%%%%%%%%%%%%%%%%%%%%%%%%%%%%%%%%%%%%%%%%%%%%%%%%%%%%%%%%%%%%%%%%%%

\subsection{Le groupo\"ide de Galois n'est pas transversalement affine}  
%-------------------------------------------------------------------%

\begin{lemme}
%-----------%
   \label{rationnaliteasl2}
%-----------%
   Soit $\widetilde{\Omega}^0$, $\widetilde{\Omega}^1$ une suite de Godbillon-Vey de type $asl_2$ pour le feuilletage donn\'e 
   par $X_1$. Il existe une suite commen\c{c}ant par $\Omega^0 = \left( \! \begin{array}{c} \omega^0_1 \\ \omega^0_2 \end{array}  \! \right)$
   pour chaque vecteur de formes v\'erifiant $\omega^0_1 \wedge \omega^0_2 =  i_{X_1} dx \wedge dy \wedge dy'$.
   S'il existe une $asl_2$-suites commen\c{c}ant par $\Omega^0$ alors elle est unique. 
%-----------%
\end{lemme} 
 
\begin{preuve}
     Les formes $\widetilde{\Omega}^0$ et $\widetilde{\Omega}^1$ satisfont les relations suivantes : 
     %aaaaaaaaaa
        $$
        \left\{
        \begin{array}{rcl}
          \widetilde{\Omega}^0(X_1) &=&0,\\
                   d\widetilde{\Omega}^0    &=&   \widetilde{\Omega}^1\wedge\widetilde{\Omega}^0,\\
                  d\widetilde{\Omega}^1    &= & \widetilde{\Omega}^1\wedge\widetilde{\Omega}^1,\\
          trace\ \widetilde{\Omega} ^1   &= & 0.
        \end{array}
        \right.
        $$
     %aaaaaaaaaa
     En \'ecrivant $\Omega^0= F \widetilde{\Omega}^0$ et $\Omega^1 = dF F^{-1} + F \widetilde{\Omega}^1 F^{-1}$, 
     on construit une nouvelle suite de Godbillon-Vey. De plus $\Omega^0_1 \wedge \Omega^0_2$ et 
     $\widetilde{\Omega}^0_1 \wedge \widetilde{\Omega}^0_2$ sont deux
     formes volumes trans\-ver\-ses in\-va\-ri\-an\-tes sous le flot des champs de vecteurs tangents. On en d\'eduit que  
     %aaaaaaaaaa
        $$
        \omega^0_1 \wedge \omega^0_2 = c\  \widetilde{\omega}^0_1 \wedge \widetilde{\omega}^0_2
        $$
     %aaaaaaaaaa 
     o\`u $c$ est une int\'egrale premi\`ere rationnelle donc une constante. Calculons la trace de $\Omega^1$: 
     %aaaaaaaaaa
        $$
        \begin{array}{rl}
          trace\ \Omega^1 &= trace\ dFF^{-1}\\
                              &= \frac{d(det F)}{detF}\\
                              &=\frac{dc}{c}=0
        \end{array}
        $$
     %aaaaaaaaaaa
     et $\Omega^0,\ \Omega^1$ forme une $asl_2$-suite.
     
     Supposons maintenant qu'il existe deux $asl_2$-suites : $\Omega^0,\ \Omega^1$ et $\Omega^0,\ \widetilde{\Omega}^1$.
     Le groupo\"ide de Galois du feuilletage est un sous-groupo\"ide de Lie du groupo\"ide d\'efinie par une de ces suites.
     Or le lemme \ref{lemmedeitchiban} montre qu'un tel feuilletage est contenu dans un feuilletage de codimension un, ce qui contredit la
     proposition \ref{prim}.       
\end{preuve}

\begin{proposition}
%------------%
   \label{asl2}
%------------%   
    Le feuilletage donn\'e par $X_1$ n'admet pas de suite de Godbillon-Vey de type $asl_2$.
%------------%
\end{proposition}

Le lemme pr\'ec\'edent affirme que si le feuilletage donn\'e par $X_1$ admet une $asl_2$-suite alors celle-ci est rationnelle.

\begin{lemme}
%------------%   
   Si elle existe, la $asl_2$-suite  commen\c{c}ant par $\widetilde{\Omega}^0 = \left( \! \! \begin{array}{c} dy'-(6y^2+x)dx \\ dy-y'dx \end{array} \! \! \right)$ 
   doit \^etre polynomiale.
%------------%
\end{lemme}

\begin{preuve}    
     Nous allons montrer que le lieu des p\^oles de la forme $\widetilde{\Omega}^1$ est une hypersurface invariante pour le feuilletage.
      Le lemme \ref{hypeinv} permettra de conclure.
     Soit $Z$ l'hypersurface des p\^oles de $\widetilde{\Omega}^1$. On construit un syst\`eme d'int\'egrales premi\`eres locales 
     (en dehors de $Z$) en r\'esolvant les syst\`emes diff\'erentiels suivant :
     %aaaaaaaaa
        $$
        \left\{
         \begin{array}{rcl}
                      dH & = & \partial H \widetilde{\Omega}^0 \\
           d\partial H & = & - \partial H \widetilde{\Omega}^1 ,
         \end{array}
        \right.
        $$
     %aaaaaaaa   
     o\`u $H = \left( \! \begin{array}{c} H_1 \\ H_2 \end{array} \! \right)$ et $\partial H$ est une matrice $2 \times 2$.
     Si $Z$ est transverse aux feuilles, les int\'egrales premi\`eres se prolongent \`a $Z$. Comme de plus 
     $dH_1 \wedge dH_2 = \left( \det \partial H \right) \widetilde{\omega}^0_1 \wedge \widetilde{\omega}^0_2$, $\det \partial H$ ne s'annule que
      le long de feuilles ;
      on en deduit que $Z$ doit \^etre tangent au feuilles. Ceci contredit le corollaire \ref{hypeinv}.    
\end{preuve}

\begin{corollaire}
\label{rationnalite}
    Si elle existe, la $asl_2$-suite commen\c{c}ant par $\Omega^0 = \left( \! \begin{array}{c} y'dy'-(6y^2+x)dy \\ \frac{dy}{y'}-dx \end{array} \! \right)$ 
    est compos\'ee de 1-formes \`a coefficients dans $\CC [x,y,y',1/y']$.
\end{corollaire}

    Soit $\Omega^0$, $\Omega^1$ cette suite. On la prolonge au param\`etre $\alpha$ de mani\`ere \`a avoir une 
    $asl_2$-suite, $\Omega^0|_\alpha$, $\Omega^1|_\alpha$ pour $X_\alpha$. 
    Ces 1-formes v\'erifient 
    %aaaaaaaaa
    \begin{eqnarray}
       d \Omega^0|_\alpha &=& \Omega^1|_\alpha \wedge \Omega^0|_\alpha \label{asl2:1}\\
       d \Omega^1|_\alpha &=& \Omega^1|_\alpha \wedge \Omega^1|_\alpha  \label{asl2:2}
    \end{eqnarray}
    %aaaaaaaaa
    On a $p_\Sigma(\Omega^0|_\alpha) = \left( \! \begin{array}{c} -6 \\ 1\end{array} \! \right)$ et
    $p_\Sigma(\Omega^1|_\alpha) = \left( \! \begin{array}{cc} 0 & -7\\ 7 & 0\end{array} \! \right)$. 
    La matrice $\Omega^1|_\alpha$ que nous obtenons ainsi a \textit{a priori} un p\^ole le long de $\alpha=0$. 
    Nous allons montrer que ce p\^ole est n\'ecessairement d'ordre z\'ero puis conclure avec le lemme \ref{asl20}.
    D\'eveloppons la matrice $\Omega^1_\alpha$ en puissance de $\alpha$ :
       $$
       \Omega^1_\alpha = \frac{1}{\alpha^n}\Xi_n + \frac{1}{\alpha^{n-1}}\Xi_{n-1}+ \ldots 
       $$ 
    
\begin{lemme}
%--------------%
\label{pas1}
%--------------%
  Si $n>0$ alors on est dans un des deux cas suivants :
  \begin{enumerate}
   \item $\Xi_{n} = \left(\begin{array}{cc} 1  &   \frac{b}{a} \\ -\frac{a}{b} & -1\end{array} \right) (a\, \omega^{0}_1|_0 + b\, \omega^{0}_2|_0)$
   avec $a = \frac{c}{9} u^{k-2}(x+2y/y')^2$ et $ b = -\frac{c}{3} u^{k-1}(x+2y/y')$ o\`u $c \in \CC$ et $k=\frac{5n+8}{6} \in \ZZ$,
   \item $\Xi_{n} = \left(\begin{array}{cc} 0 & 0 \\ cu^k & 0 \end{array}\right)\omega^{0}_1|_0$ avec $k=\frac{5n - 13}{6} \in \ZZ$.
  \end{enumerate}
%--------------%
\end{lemme}
\begin{remarque}
  Si $n>0$ alors $n>1$. Le seul cas posant probl\`eme est le cas \emph{(2)}, $n=1$ et $k=-2$. En posant $\beta=\alpha^5$, on obtient un champs 
  $X_{\beta^{1/5}}$ donnant $X_1$ en $\beta=1$ et une $asl_2$-suite donn\'ee par $\Omega^1_\beta = \frac{1}{\beta^{n/5}}\Xi_n + \ldots $.
  Le corollaire \ref{rationnalite} implique que $n$ est divisible par $5$.  
\end{remarque}

\begin{preuve}
    Posons $\Xi_{n} = \left(\begin{array}{cc} \eta_{11} & \eta_{12} \\ \eta_{21} & -\eta_{11} \end{array}\right)$.
    Nous allons distinguer deux cas suivant que $\eta_{11}$ soit nulle ou non.\\
    
    \noindent
    \underline{Premier cas, $\eta_{11} \not = 0$.}\\
    %------------------------------------------------%
    
    L'\'egalit\'e $\Xi_n \wedge \Xi_n$ donne $\Xi_{n} = \left(\begin{array}{cc} 1& f \\ g & -1 \end{array}\right)\eta_{11}.$
    Comme de plus $\Xi_n \wedge \Omega^0|_0 = 0$, 
    %aaaaaaaaaaa
       $$
       \Xi_{n} = \left(\begin{array}{cc} 1& b/a \\ -a/b & -1 \end{array}\right)(a\, \omega_1^0|_0 + b\, \omega_2^0|_0).
       $$
    %aaaaaaaaaaa
    Consid\'erons ensuite l'\'egalit\'e    
    %aaaaaaaaaa
     \begin{equation}
       \label{equadiff}
        i_{X_{\alpha}}d\Omega^{1}|_\alpha = \left[ \Omega^{1}|_{\alpha}(X_\alpha) , \Omega^{1}|_{\alpha} \right] 
     \end{equation}
    %aaaaaaaaaa
    avec $\Omega^{1}|_{\alpha}(X_\alpha)= \left(\begin{array}{cc} 0 & -\alpha y' \\ \frac{1}{y'^2} & 0 \end{array}\right).$
    En calculant le coefficient de $1 / \alpha^{n}$, on obtient
    %aaaaaaaaaa
       $$
       i_{X_0}d\,\Xi_{n} = \frac{1}{y'^2}\left[\left(\begin{array}{cc} 0 & 0 \\ 1 & 0 \end{array}\right), \Xi_{n}\right]
       $$
    %aaaaaaaaaa
    Calculons les deux membres de l'\'egalit\'e \\
    %aaaaaaaaaa
       $
       \begin{array}{lcl}
       \bullet \ \ i_{X_0}d\,\Xi_{-n} &=& 
          \left( \begin{array}{cc} 0 & X_0 \frac{b}{a}\\ -X_0 \frac{a}{b} & 0 \end{array}\right) (a\, \omega^{0}_1|_0 + b\, \omega^{0}_2|_0)
       + \left(\begin{array}{cc} 1& \frac{b}{a} \\ -\frac{a}{b} & -1 \end{array}\right) ((X_0a + \frac{b}{y'^2}) \omega^{0}_1|_0 + X_0b\ \omega^{0}_2|_0),
       \end{array}
       $\\
       $
       \begin{array}{lcl}
         \bullet \ \  \left[\left(\begin{array}{cc} 0 & 0 \\ 1 & 0 \end{array}\right), \Xi_{n}\right]
                      &=& \left(\begin{array}{cc} -\frac{b}{a} & 0 \\ 2& \frac{b}{a}  \end{array}\right) (a\, \omega^0_1|_0 +b\, \omega^0_2|_0).
       \end{array}
       $\\
    %aaaaaaaaaa
    En identifiant les coefficients de $\omega^{0}_1|_0$ et $\omega^{0}_2|_0$ de ces matrices, on obtient le syst\`eme
    %aaaaaaaaaa
       $$
       \left \{
        \begin{array}{l}
           X_0 a = \frac{-2b}{y'^2}\\
           X_0 b = -\frac{b^2}{a} \frac{1}{y'^2}.
        \end{array}
       \right.
       $$
    %aaaaaaaaaa
    On en d\'eduit que $\frac{b^2}{a}$ est une int\'egrale premi\`ere homog\`ene de $X_0$
    donc $\frac{b^2}{a} = c u^k$ avec $c \in \CC$ et $k \in \ZZ$. On r\'esoud le 
    syst\`eme en utilisant les formules (\ref{formules}).
    En calculant le poids de $\Xi_n$ \`a partir des formules d'un cot\'e et \`a partir du poids de $\Omega_\alpha$ 
    de l'autre, on obtient la condition sur l'ordre du p\^ole en $\alpha$.\\
    
    \noindent
    \underline{Deuxi\`eme cas, $\eta_{11}=0$.}\\
    %--------------------------------------------%
    
    Comme $\Xi_n \wedge \Xi_n = 0$ et $\Xi_n \wedge \Omega^{0}|_0 = 0$, on a soit 
    $\Xi_n = \left(\begin{array}{cc} 0 & 0 \\ g & 0 \end{array}\right)\omega_1^{0}|_0 $, soit 
    $\Xi_n = \left(\begin{array}{cc} 0 & f \\ 0 & 0 \end{array}\right)\omega_2^0|_0 $.
    En reprenant les \'equations diff\'erentielles donn\'ees par (\ref{equadiff}), on obtient le r\'esultat annonc\'e. 
\end{preuve}
    
\begin{preuved}{de la proposition \ref{asl2}}
 Nous allons maintenant calculer $\Xi_{n-1}$ et aboutir \`a une contradiction.\\
    
    \noindent
    \underline{Premier cas, $\eta_{11} \not = 0$.}\\
    %----------------------------------------------%
    %
    
    Dans ce cas, $\Xi_n$ est \'egale \`a $\left( \! \begin{array}{cc} 1& b/a \\ -a/b & -1 \end{array} \! \right) \eta_{11}$.
    En \'ecrivant le coefficient de $\frac{1}{\alpha^{n-1}}$ dans (\ref{asl2:1}), on obtient
    %aaaaaaaaaaa
       $$
       \Xi_{n-1} \wedge \Omega^0|_0 + \Xi_{n} \wedge \left( \! \begin{array}{c} -xdy \\ 0 \end{array} \! \right) = 0.
       $$
    %aaaaaaaaaa
    Ce qui implique que 
    %aaaaaaaaaa
       $$
       \Xi_{n-1} = A\, \omega^0_1|_0 + B\, \eta_{11} - \left( \! \begin{array}{cc}    
                                                                                                a        &  b \\
                                                                                    -\frac{a^2}{b} & -a  
                                                                                   \end {array}     \!        \right) xdy
       $$
    %aaaaaaaaaa
    avec $A \left( \! \begin{array}{c} 0 \\ 1 \end{array} \! \right) + B \left( \! \begin{array}{c} -b \\ a \! \end{array} \right) = 0$.
    %aaaaaaaaaaa
    En \'ecrivant le coefficient de $\frac{1}{\alpha^{2n-1}}$ dans (\ref{asl2:2}), on obtient 
    %aaaaaaaaaaa
       $$
       \Xi_n \wedge \Xi_{n-1} + \Xi_{n-1} \wedge \Xi_n = 0
       $$
    %aaaaaaaaaa
    Ce qui donne $\left[ A, \left(\begin{array}{cc} 1& b/a \\ -a/b & -1 \end{array}\right) \right] = 0$, c'est-\`a-dire
       $$
       A = \ell \left(\begin{array}{cc} 1& b/a \\ -a/b & -1 \end{array}\right)
       $$
    avec $\ell$ quelconque.
    On en d\'eduit ensuite que
       $$
       B= \ell \left(\begin{array}{cc} 0& - b/a^2 \\ -1/b & 0 \end{array}\right) 
          + t \left(\begin{array}{cc} 1& b/a \\ -a/b & -1 \end{array}\right)
       $$
    avec $t$ quelconque.
    Consid\'erons maintenant l'\'equation diff\'erentielle portant sur $\Xi_{n-1}$ que l'on d\'eduit de (\ref{equadiff}) :
    %aaaaaaaaaa
       $$
         i_{X_0}d\,\Xi_{n-1} + i_{x\frac{\partial}{\partial y'}} d\,\Xi_n
       = \frac{1}{y'^2}\left[\left(\begin{array}{cc} 0 & 0 \\ 1 & 0 \end{array}\right), \Xi_{n-1}\right]
       - y' \left[\left(\begin{array}{cc} 0 & 1 \\ 0 & 0 \end{array}\right), \Xi_{n}\right].
       $$
    %aaaaaaaaaa
    Calculons les termes ci-dessus :\\
    
    \noindent
    $
    \begin{array}{rcl}
     \bullet \ i_{X_0}d\,\Xi_{n-1}  & = &  
            \left(   
            \begin{array}{cc} 
          X_0 \ell + a X_0 t                                                 &  -\left(\frac{b}{a}\right)^2 \frac{\ell}{y'^2} + b X_0 t \\
        -2\frac{a}{b} X_0 \lambda + \frac{\ell}{y'^2} - \frac{a^2}{b} X_0 t &  -X_0 \ell - a X_0 t
             \end{array} 
             \right)    
              \omega_1^0|_0  \\
            & & \\
        &+& \left( 
            \begin{array}{cc} 
                     bX_0 t -ay'                     & -\left(\frac{b}{a}\right)^2 X_0 \ell -\ell \left(\frac{b}{a}\right)^3\frac{1}{y'^2}+ \frac{b^2}{a} X_0 t -by'\\
        -X_0 \ell -a X_0 t + \frac{a^2}{b} y' & -b X_0 t +ay'
             \end{array} 
             \right) 
             \omega_2^0|_0 \\
             & & \\
        &+& \left( 
             \begin{array}{cc}
                   2b  &  \frac{b^2}{a} \\
                  -3a  &  -2b     
             \end{array} 
             \right)
              \frac{x}{y'^2}dy 
           + xy' \left( 
                \begin{array}{cc}
                        da    &  db \\
                 -d\left(\frac{a^2}{b}\right) & -da           
                \end{array} 
                \right) 
           + t \ i_{X_0}d\Xi_n ,
      
    \end{array} 
    $\\
     
    %aaaaaaaaaa
    %aaaaaaaaaa
    \vspace{0.6cm}
    \noindent
       $
       \begin{array}{lcl}
      \bullet \ i_{x\frac{\partial}{\partial y'}} d\,\Xi_n & = &

        x \left(   
                   \begin{array}{cc} 
                      \frac{\partial a}{\partial y'}        &  \frac{\partial b}{\partial y'}  \\
                   - \frac{\partial (a^2/b)}{\partial y'}  &  -\frac{\partial a}{\partial y'} 
                   \end{array} 
                   \right) 
                   \omega_1^0|_0 
           + x \left( 
                \begin{array}{cc} 
                    \frac{\partial b}{\partial y'}   &  \frac{\partial (b^2/a)}{\partial y'}  \\
                  - \frac{\partial a}{\partial y'}   &  - \frac{\partial b}{\partial y'} 
                \end{array}  
                \right) 
                \omega_2^0|_0 \\
                & & \\
             
             &+& \left( 
                 \begin{array}{cc}
                    -b  &  -\frac{b^2}{a} \\
                     a  &   b     
                 \end{array} 
                 \right) 
                 \frac{x}{y'^2}dy 
           -xy' \left( 
                 \begin{array}{cc}
                         da    &  db \\
                  -d\left(\frac{a^2}{b}\right) & -da           
                \end{array} 
                \right),
       \end{array}
       $\\
       
    %aaaaaaaaaa
       \vspace{0.6cm}
       \noindent
       $
       \begin{array}{rcl}
       \bullet \ \left[\left(\begin{array}{cc} 0 & 0 \\ 1 & 0 \end{array}\right), \Xi_{n-1}\right] &=&
         \ell \left(  
                    \begin{array}{cc} 
                        0    &  0  \\
                        2    &  0 
                    \end{array} 
                    \right) 
                    \omega_1^0|_0 
       + \ell \left( 
                     \begin{array}{cc} 
                        \left(\frac{b}{a}\right)^2   &  0  \\
                           0         &  - \left(\frac{b}{a}\right)^2 
                     \end{array} 
                     \right) 
                     \omega_2^0|_0 
         +         \left( 
                   \begin{array}{cc}
                       b  &  0 \\
                    -2a  & -b     
                  \end{array} 
                  \right) 
                  x dy 
         + t  \left[ 
                  \left( 
                  \begin{array}{cc}
                         0    &  0 \\
                         1 & 0           
                  \end{array}
                  \right)
                   ,  \Xi_n \right]
       \end{array}                     
       $\\
       
    %aaaaaaaaaa
       \vspace{0.6cm}
       \noindent
       $
       \begin{array}{l}
       \bullet \ \left[\left(\begin{array}{cc} 0 & 1 \\ 0 & 0 \end{array}\right), \Xi_{n}\right]=
        \left( 
       \begin{array}{cc}
              -a/b  &  -2 \\
                0    &  a/b     
        \end{array} 
        \right) 
        (a\omega_1^0|_0 + b \omega_2^0|_0). 
        \end{array}
        $ \\
    
    %aaaaaaaaaa
    
    Parmis les \'equation que l'on obtient en identifiant les coefficients, on trouve les deux \'equations suivantes 
    %aaaaaaaaaa
       $$
       \left\{
       \begin{array}{l}
         X_0 t = 2\frac{a}{b}y' + \frac{b}{a^2} \frac{\ell}{y'^2} - \frac{\partial b}{ \partial y'}\frac{x}{b}\\
         X_0 \ell = - \frac{b}{a}\frac{\ell}{y'^2} - \frac{a^2}{b}y' - b \frac{\partial}{\partial y'}\left(\frac{a}{b}\right) x 
       \end{array}
       \right.
       $$
    %aaaaaaaaaa
    En rempla\c{c}ant $a$ et $b$ par leurs valeurs (lemme \ref{pas1}), on obtient l'\'equation suivante pour $\ell$ :
    $$
    X_0\ell = \frac{3u}{x+2\frac{y}{y'}} \frac{\ell}{y'^2} + \frac{c}{27}u^{k-3} \left( x+2\frac{y}{y'} \right)^3 y'
                     + \frac{2c}{9} u^{k-3}x\left( x+2\frac{y}{y'} \right)^2 y'+ \frac{2c}{9} u^{k-2}x\left( x+2\frac{y}{y'} \right) \frac{y}{y'^2}.
    $$
    On r\'esoud cette derni\`ere \'equation en posant $\ell = \lambda \frac{x+2\frac{y}{y'}}{3u} \frac{9}{c u^{k-2}}$. On a alors
    $$
    X_0 \lambda = 7 y' x^2 + \left( 16 y + 6 u \frac{y}{y'^2}\right) x + \frac{4y^2}{y'}.
    $$ 
    La fonction $\lambda$ est donc un polyn\^ome en $x$ de degr\'e inf\'erieur \`a trois. Par homog\'en\'eit\'e sous $\Sigma$, $\lambda$ est de degr\'e $2$. 
    On r\'esoud les \'equations portants sur les coefficients de $\lambda$ et on obtient finalement
    $$
    \ell = \frac{x+2\frac{y}{y'}}{3u} \frac{9}{c u^{k-2}} \left( 7yx^2 + 4\frac{y^2}{y'} x +d\right) 
    $$
    o\`u $d$ est constante. 
    L'\'equation satisfaite par $t$ est alors
    \begin{multline}
    X_0t  = - \left(\frac{y^3}{y'^2} + d \right) \frac{1}{y'^2} \frac{1}{\left( x+2\frac{y}{y'}\right)^2}  
           +  \left( 28 \frac{y^2}{y'}\right)\frac{1}{y'^2}\frac{1}{x+2\frac{y}{y'}} \nonumber \\
           + 4(k-1)\frac{y}{u} - 9\frac{y}{y'^2} - (2k-\frac{8}{3})\frac{y'}{u}\left( x+2\frac{y}{y'}\right).
    \end{multline}  
    Puisque $X_0 t$ a un p\^ole d'ordre $2$ le long de $x+2\frac{y}{y'} = 0$, $t$ doit avoir un p\^ole d'ordre $1$ :
    $$
    t = \frac{t_{-1}}{x+2\frac{y}{y'}}+t_0+t_1\left( x+2\frac{y}{y'}\right)+ t_2 \left( x+2\frac{y}{y'}\right)^2
    $$
    On a alors $t_{-1} = -\frac{1}{3u}\left( 48 \frac{y^3}{y'} + d \right)$ et $X_0 t_{-1} = \frac{28}{y'^2}\frac{y^2}{y'}$, ce qui est impossible. Il n'y a pas
    de solution polynomiale \`a l'\'equation portant sur $t$. Il n'existe donc pas de $asl_2$-suite pour le feuilletage avec $\eta_{11} \not = 0$. \\
    
    \noindent
    \underline{Deuxi\`eme cas, $\eta_{11}$ est nulle.}\\
    %----------------------------------------------------%
    
    Dans ce cas, $\Xi_n = \left( \begin{array}{cc} 0 & 0 \\ cu^k & 0  \end{array}\right) \omega_1^0|_0$. Calculons $\Xi_{n-1}$.
    Les conditions $\Xi_n \wedge \Xi_{n-1} + \Xi_{n-1} \wedge \Xi_{n} = 0$ et 
    $\Xi_{n-1} \wedge \omega^0|_0 + \Xi_{n} \wedge \left( \begin{array}{c} -xdy \\ 0 \end{array} \right)= 0$
    donnent 
    %aaaaaaaaaa
       $$
       \Xi_{n-1} = \left(\begin{array}{cc} -e & 0 \\ d & e \end{array}\right)\omega_1^0|_0 
       + \left(\begin{array}{cc} 0 & 0 \\ e & 0 \end{array}\right)\omega_2^1|_0 
       - \left(\begin{array}{cc} 0 & 0 \\ cu^k & 0 \end{array}\right)xdy.
       $$
    %aaaaaaaaaa
    En reprenant l'\'equation (\ref{equadiff}) et en remarquant que $\Xi_{n}$ est ferm\'ee, on obtient
    %aaaaaaaaaa
       $$
       i_{X_0}d\,\Xi_{n-1} = \frac{1}{y'^2} \left[\left(\begin{array}{cc} 0 & 0 \\ 1 & 0  \end{array}\right) , \Xi_{n-1}  \right] - y'\left[\left(\begin{array}{cc} 0 & 1 \\ 0 & 0  \end{array}\right) , \Xi_{n}  \right]
       $$
       qui donne
       \begin{multline}
      \left(\begin{array}{cc} -X_0e & 0 \\ X_0d & X_0e  \end{array}\right)\omega_1^0|_0 + \left(\begin{array}{cc} 0 & 0 \\ X_0e & 0  \end{array}\right)\omega_2^0|_0 + \left(\begin{array}{cc} 0 & 0 \\ e & 0  \end{array}\right)\frac{1}{y'^2}\omega_1^0|_0 - \left(\begin{array}{cc} 0 & 0 \\ cu^k & 0  \end{array}\right)y' \omega_2^0|_0 \nonumber \\      
      = y'\left(\begin{array}{cc} -cu^k & 0 \\ 0 & cu^k \end{array}\right)\omega_1^0|_0 
                              + \frac{1}{y'^2}\left(\begin{array}{cc} 0 & 0 \\ 2e & 0 \end{array}\right)\omega_1^0|_0.
       \end{multline}
    %aaaaaaaaaa
    Ce qui donne les \'equations suivantes $X_0e = cu^k y'$ et $X_0 d = \frac{3e}{y'^2}$, c'est-\`a-dire $c= cu^ky$ 
    et $X_0 d = 3cu^k\frac{y}{y'^2}$.
    %
    %------------%
    \begin{lemme}
    %------------%
        L'\'equation $X_0 R = \frac{y}{y'^2}$ n'a pas de solution rationnelle.
    %------------%    
    \end{lemme}
    %------------%
    %
    \begin{preuve}
    %------------%
          On d\'eduit de cette \'equation que $X_0 \frac{\partial R}{\partial x} = 0$. D'un autre cot\'e le poids sous $\Sigma$
          de $R$ doit \'egale \`a 4 donc  $\frac{\partial R}{\partial x} = 0$. En \'ecrivant l'\'equation dans les coordonn\'ees
          $x,y,u$, on obtient :
          %aaaaaaaaaa
             $$
             (4y^3 + u)\frac{\partial R_1}{\partial y} + 6 y^2 R_1 = \frac{y}{4y^3 + u},
             $$
          %aaaaaaaaaa
          o\`u $R=R_0(y,u) + y'R_1(y,u)$.
          Si cette \'equation admet une solution rationnelle, elle s'\'ecrit $R_1=\frac{P}{4y^3+u}$ avec $P$ polynomiale. 
          En calculant le coefficient du terme de plus haut degr\'e en $y$ dans l'\'equation satisfaite par $P$,
          on v\'erifie $P$ ne peux pas \^etre polynomiale. 
    %-----------%
    \end{preuve}
    
    L'\'equation sur $d$ n'a de solution rationnelle que si $c$ est nul, ce qui contredit le choix de $n$. La matrice de 1-forme
    $\Omega^1|_{\alpha}$ ne peux donc pas avoir de p\^ole en $\alpha$. Mais d'apr\`es le lemme \ref{asl20}, elle ne peux pas ne pas 
    avoir de p\^ole. Il n'existe donc pas de matrice de formes \`a coefficients rationnels satisfaisant les \'equations de
    $\Omega^1|_{\alpha}$. Ceci prouve qu'il n'existe pas de $asl_2$-suite pour le feuilletage donn\'e par la premi\`ere \'equation 
    de Painlev\'e.  
%-------------------------------------------------------------------% 
\end{preuved}
%-------------------------------------------------------------------%

\section{Irr\'eductibilit\'e de $P_1$}

Dans cette section, nous allons utiliser le groupo\"ide de Galois du feuilletage d\'efini par $P_1$ pour montrer son ir\-r\'e\-duc\-ti\-bi\-li\-t\'e. 
Nous commencerons par d\'efinir ce que nous entendons par feuilletage r\'eductible. En utilisant les notions de ``type diff\'erentiel'' et de ``degr\'e 
typique de transcendance diff\'erentielle'' de Kolchin, nous montrerons ensuite que le groupo\"ide 
de Galois d'un feuilletage r\'eductible est plus petit que le groupo\"ide de Galois de $P_1$.

\subsection{Feuilletages r\'eductibles}

Dans \cite{Umemura-painleve}, une solution particuli\`ere d'une \'equation dif\-f\'e\-ren\-tiel\-le d'ordre deux est dite r\'eductible si on peut l'exprimer 
rationnellement apr\`es avoir r\'esolu successivement des \'equations diff\'erentielles lin\'eaires (ou associ\'ees \`a un groupe alg\'ebrique) 
et des \'equations d'ordre un. Cette d\'efinition d'ir\-r\'e\-duc\-ti\-bi\-li\-t\'e ne concerne que les solutions particuli\`eres de l'\'equation ind\'ependemment
de l'\'equation.  Dans l'esprit de la d\'efinition du groupo\"ide de Galois, nous allons d\'efinir une notion de r\'eductibilit\'e 
``globale" du feuilletage \textit{i.e.} de r\'eductibilit\'e de la solution g\'en\'erale. 

Soit $\F$ un feuilletage de codimension deux. Les feuilletages que nous consid\`ererons comme plus simples que $\F$ sont 
les feuilletages de codimension un d'une part et les feuilletages donn\'es par une connexion lin\'eaire (ou associ\'ee \`a un groupe alg\'ebrique) de 
l'autre ainsi que leurs versions relative le long d'un feuilletage de codimension un. Dans un autre contexte, ce type d'extension est \'etudi\'e dans 
\cite{cassidy-singer}.

Soient $(K,\partial_1,\ldots,\partial_n)$ un corps diff\'erentiel. On notera $T_K$ le $K$-espace vectoriel engendr\'e par les $\partial$, $T^*_K$ son 
dual et $\{d_1, \ldots, d_n\}$ la base duale des $\partial$ 

\begin{definition}[\cite{cassidy-singer}]
\label{modulo}On dira que une extension diff\'erentielle $K \subset L$ est une extension 
fortement normale de $K$ relative en codimension $q$ si il existe un sous-corps $Q$ de $K$ engendr\'e par $q$ \'el\'ements fonctionnellement 
ind\'ependants tel que si on note $\underline{K}$ le corps $K$ muni des d\'erivations $T_{\underline{K}} = \cap_{\atop h \in Q} \ker dh$ et 
$\overline{Q}^{alg}$ la cl\^oture alg\'ebrique de $Q$, l'extension 
$\underline{K} \otimes \overline{Q}^{alg} \subset \underline{L} \otimes \overline{Q}^{alg}$ est fortement normale.
\end{definition}

Une extension fortement normale \`a un groupe de Galois. Dans \cite{cassidy-singer} les derivations n'appartenant pas \`a  $T_{\underline{K}}$ 
sont prise en consid\'eration pour \'etudier l'extension ce qui permet de d\'efinir un groupe de Galois pour l'extension $K \subset L$, c'est un groupe 
alg\'ebrique diff\'erentiel.   

\begin{exemple}
Consid\'erons le champs de vecteurs $X_0 = \frac{\partial}{\partial x} + y'\frac{\partial}{\partial y} + 6y^2\frac{\partial}{\partial y'}$. Ce champs
admet une int\'egrale premi\`ere rationnelle $u = y'^2 - 4y^3$ et une int\'egrale premi\`ere dans une extension fortement normale relative en 
codimension un de $\CC(x,y,y')$ : $ H = x - \int \frac{dy}{\sqrt{4y^3 + u}}$. Son groupe de Galois est un sous-groupe alg\'ebrique diff\'erentiel de 
$\GG_a(\overline{\CC(u)}^{alg})$.
\end{exemple}

\begin{definition}
\label{reductibilite}
   Un feuilletage $\F$ de codimension deux sur $\CC^n$ est dit r\'eductible 
   si il existe deux integrales premi\`eres $H$ dans une extension diff\'erentielle $K_p$ de $\CC(x_1,\ldots,x_n)$ construite de la mani\`ere suivante
     $$
      \CC(x_1,\ldots,x_n) = K_0 \subset K_1 \ldots \subset K_p
     $$    
   telle que les extensions interm\'ediaires $K_i \subset K_{i+1}$ sont
   \begin{itemize}
      \item des extensions alg\'ebriques, 
      \item des extension fortement normales \cite{Kolchin} 
      \item des extensions fortement normales relatives en codimension un. 
       \item des extensions par une int\'egrale
      premi\`ere non constante d'un feuilletage de codimension un, \textit{i.e.} $K_{i+1} = K_i (<H>)$ avec $dH \wedge \omega =0$ pour une 1-forme int\'egrable
      \`a coefficients dans $K_i$. 
     
    \end{itemize}
\end{definition}   
Cette d\'efinition se g\'en\'eralise aisement aux feuilletages de codimension quelconque.
Dans le cadre de cet article, le th\'eor\`eme que nous allons prouv\'e est le suivant.
\begin{theoreme}
\label{irreductibilite}
   Le feuilletage d\'efini par $P_1$ est irr\'eductible.
\end{theoreme} 

La preuve est donn\'ee dans la section suivante. Dans un premier temps nous prouvons que les feuilletages reductibles admettent des int\'egrales premieres
satisfaisant un ''gros'' syst\`eme d'\'equations aux d\'eriv\'ees partielles. La taille de ``l'espace des solutions'' de ce syst\`eme est mesur\'e par le 
``type d'une extension diff\'erentielle'' d\'efini dans la section suivante. Le pseudo-groupe d'holonomie du feuilletage fixant les int\'egrales 
premi\`eres locales (lorsque cela \`a un sens), ces \'elements satisfont aussi un syt\`eme d'e.d.p. La taille du groupo\"ide de Galois donne une 
borne inf\'erieure \`a la taille de ce syst\`eme d'e.d.p. Dans le cas du feuilletage donn\'e par $P_1$, cette borne est sup\'erieure \`a la taille de 
l'espace des int\'egrales premi\`eres particuli\`eres d'un feuilletage r\'eductible.

\subsection{Types d'une extension diff\'erentielle et preuve du th\'eor\`eme \ref{irreductibilite}}

La d\'efinition du type d'une extension de corps diff\'erentiels est bas\'ee sur un analogue diff\'erentiel du polyn\^ome de Hilbert pour les 
$\D$-vari\'et\'e introduit par Kolchin sous le nom de ``polyn\^ome de transcendance diff\'erentielle'' (\cite{Kolchin}, chapitre II).
Rappellons d'abord la d\'efinition du type d'une $\D$-vari\'et\'e de sections de $Z$ sur $X$. 
\begin{definition}
   Soient $Z$ une vari\'et\'e sur $X$ et $\Y$ une $\D$-vari\'et\'e irr\'eductible de $J(Z_{/X})$. On d\'efinit la croissance de $\Y$ comme 
   la suite d'entier $c_\ell = \dim_X \Y_\ell$. Le type de $\Y$ est le monome $a\ell^b$ tel que $c_\ell \sim_{\infty} a\ell^b $.
\end{definition}

Dans \cite{Kolchin} chapitre II-12, le polyn\^ome interpolant la suite $c_\ell$ pour de grandes valeurs de $\ell$ est appell\'e ``polyn\^ome de 
transcendance diff\'erentielle''. Dans \cite{Kolchin} chapitre II-13, les entiers $a$ et $b$ sont appell\'es ``degr\'e typique de transcendance 
diff\'erentielle'' et ``type diff\'erentiel".    

\begin{exemple}
   La croissance de la $\D$-vari\'et\'e $J(\CC^n \to \CC^p )$ est $c_\ell = p\frac{(n+\ell)!}{n! \ell !}$. 
   Son type est $\frac{p}{n!} \ell ^n$.
\end{exemple}
\begin{exemple}
   Soient $\F$ un feuilletage de codimension un de $\CC^n$ donn\'e par une forme $\omega$ et $IP(\F)$ la $\D$-vari\'et\'e des int\'egrales 
   premi\`eres \textit{i.e.} la sous-$\D$-vari\'et\'e de $J(\CC^n \to \CC)$ d\'efinie par l'id\'eal diff\'erentiel engendr\'e par les composantes de
   $dH \wedge \omega$. La croissance de $IP(\F)$ est $c_\ell = \ell +1$. 
\end{exemple}   
\begin{exemple}
   Soit $\F$ un feuilletage de codimension deux de $\CC^n$. La croissance de la $\D$-vari\'et\'e des couples int\'egrales premi\`eres est 
   $c_\ell = (\ell +2)(\ell +1)$. Si de plus le feuilletage pr\'eserve une forme volume transverse, la croissance de la $\D$-vari\'et\'e 
   des couples d'int\'egrales premi\`eres compatibles au volume est $c_\ell =\frac{1}{2} (\ell+2)(\ell+1)+\ell +1$.   
\end{exemple}
\begin{exemple}
   Soit $V$ un espace vectoriel muni d'une connexion int\'egrable sur $X$. La croissance de la $\D$-vari\'et\'e des sections plates de $V$ est 
   $c_\ell = \dim_X V$. Elle est ind\'ependante de $\ell$.
\end{exemple}
 
La d\'efinition donn\'ee par Kolchin dans le cadre des extensions de corps diff\'erentiels se d\'eduit de la d\'efinition p\'rec\'edente de la mani\`ere 
suivante. 
Soient $K, (\partial_1,\ldots, \partial_n)$ un corps diff\'erentiel et $L$ une extension diff\'erentiel de $K$. Choisissons $p$ \'el\'ements 
$y_1,\ldots,y_p$ engendrant diff\'erentiellement $L$ sur $K$. Soit $\I$ le noyau du morphisme $K<Y_1,\ldots,Y_p> \to L$. Le corps $L$ 
est alors le corps des fractions de la $\D$-vari\'et\'e au-dessus de $K$ d\'efinie par $\I$. Notons $\Y$ cette $\D$-vari\'et\'e.

\begin{lemme-def}[\cite{Kolchin}]
\label{type}
   Le type $t_{L/K}$ de $K \subset L$ est le type de $\Y$. Il ne d\'epend que de l'extension $K \subset L$. 
\end{lemme-def}

\begin{preuve}
   Soient $y_1,\ldots,y_p$ et $z_1, \ldots, z_q$ deux syst\`emes g\'en\'erateurs de $L$ sur $K$. Notons $\Y$ et $\Z$ les $\D$-vari\'et\'es au-dessus
   de $K$ qu'ils d\'efinissent. Par d\'efinition, il existe des fractions diff\'erentielles telles que $y_j = y_j(z)$. Soit $i$ l'ordre maximal de ces fractions.
   Elles d\'efinissent une application rationnelle dominante de $\Y_i$ sur $\Z_0$. Par d\'erivations on obtient des applications dominantes de 
   $\Y_{i+\ell}$ sur $\Z_\ell$. Ce qui implique $ c(\Y)_{i+\ell} \geq c(\Z)_{\ell}$. De la m\^eme mani\`ere, $c(\Z)_{j+\ell} \geq c(\Y)_{\ell}$. 
   Ceci prouve le lemme. 
\end{preuve}

\begin{exemple}
Le type d'une extension fortement normale est fini.
\end{exemple}

\begin{exemple}
Soient $(K, \partial_1, \ldots, \partial_n)$, $\omega$ une 1-forme sur $K$ (\textit{i.e.} un \'el\'ement du dual du $K$-espace vectoriel engendr\'e par les d\'erivations)
et $L = K(<H>)$ une extension diff\'erentielle engendr\'ee par une ind\'etermin\'ee $H$ satisfaisant $dH \wedge \omega = 0$. Le type $t_{L/K}$ est lin\'eaire.  
\end{exemple}

\begin{exemple}
   Soient $K \subset L$ une extension fortement normale en codimension un et $\partial_1,\ldots,\partial_n$ une base des derivations telle que 
   $T_{\underline{K}} = K\partial_1+\ldots+K\partial_{n-1}$. L'extension \'etant fortement normale par rapport \`a $T_{\underline{K}}$, on peut
   trouver une base $\{H_1,\ldots,H_p\}$ de $L$ sur $K$ telle que $\partial_j H_i$ soit alg\'ebrique sur $K(H_1,\ldots,H_p)$ pour tout $i$ et 
   $1 \leq j \leq n-1$. Le type de $L/K$ est donc lin\'eaire.
\end{exemple}

De la m\^eme mani\`ere que pour le lemme \ref{type}, on prouve le lemme suivant.

\begin{lemme}
\label{extension}
   Soit $K \subset K_1 \subset K_2$ une suite d'extensions diff\'erentielles. On a l'in\'egalit\'e suivante :
   $$
   t_{K_2/K} \leq_\infty t_{K_2 / K_1} + t_{K_1/K},
   $$
   o\`u $\leq_\infty$ signifie asymptotiquement inf\'erieur.
\end{lemme}
\begin{preuvedt}{irreductibilite}
%-----------------------------%
%
   Consid\'erons l'espace $J(\CC^3 \to \CC^2)$. Le $\D$-groupo\"ide de Lie $J^*(\CC^3 \to \CC^3)$ agit sur cet espace par composition \`a la source.
   Soit $\I$ un id\'eal diff\'erentiel r\'eduit de $\O_{J(\CC^3 \to \CC^2)}$ d\'efinissant une sous-$\D$-vari\'et\'e $\Y$ de cet espace. Nous appellerons
   stabilisateur de $\I$ (ou de $\Y$) et noterons $Stab(\I)$ (ou $Stab(\Y)$) le $\D$-groupo\"ide de Lie maximal laissant $\Y$ invariante. Notons $c_{source}$ 
   la composition \`a la source : 
   %aaaaaaaa
   $$
      c_{source} : \left(J^*(\CC^3 \to \CC^3),t \right) \times_{\CC^3} \left( J(\CC^3 \to \CC^2), s \right) \to J(\CC^3 \to \CC^2),
   $$
   %aaaaaaaa
   l'id\'eal du stabilisateur de $\I$ est le plus petit id\'eal diff\'erentiel $\J \subset \O_{J^*(\CC^3 \to \CC^3)}$ tel qu'en dehors d'une hypersurface
   $S \subset \CC^3$ on ait les inclusions suivantes 
   %aaaaaaa
   $$
      c_{source}^*\I + \J \subset \I + \J \textrm{  et  } c_{source}^*\I + i^*\J \subset \I + \J.
   $$
   %aaaaaaa
   \`A partir de ces inclusions, on montre que $J$ est l'id\'eal d'un $\D$-groupo\"ide de Lie (voir par exemple \cite{these}).
   Consid\'erons $IP(\F)$ la $\D$-vari\'et\'e des int\'egrales premi\`eres du feuilletage $\F$ d\'efini par $P_1$. C'est la sous-$\D$-vari\'et\'e de 
   $J(\CC^3 \to \CC^2)$ d\'ecrivant les couples d'int\'egrales premi\`eres. Elle est donn\'ee par l'id\'eal diff\'erentiel r\'eduit engendr\'e par
   les deux \'equations $X_1H_1$ et $X_1H_2$ o\`u $(H_1,H_2)$ sont des coordonn\'ees sur l'espace but. Le stabilisateur de $IP(\F)$ est le 
   $\D$-groupo\"ide de Lie $Aut(\F)$ dont les solutions sont les diff\'eomorphismes formelles $\Gamma$ satisfaisant les \'equations diff\'erentielles 
   $\Gamma^*X_1 \wedge X_1 = 0$.
   Soit $\I$ un id\'eal diff\'erentiel premier de $\O_{IP(\F)}$ et $\widetilde{\I}$ sa pr\'eimage dans $\O_{J(\CC^3 \to \CC^2)}$. Le stabilisateur de $\I$
   est par d\'efinition le $\D$-groupo\"ide de Lie $Stab(\widetilde{\I}) \cap Aut(\F)$.
   Pla\c{c}ons nous au voisinage d'un point r\'egulier de $X_1$ et choisissons une coordonn\'ee tangente $z$
   et deux coordonn\'ees transverses $t_1$ et $t_2$. Dans ces coordonn\'ees analytiques, l'id\'eal $\I$ admet un syst\`eme de g\'en\'erateurs ne 
   d\'ependant 
   pas de $z$ (\cite{these,codimun}). Tous les flots locaux de $X_1$ laissent donc $\I$ invariant. Autrement dit, $Stab(\I)$ est localement admissible, 
   donc admissible.   
   On obtient l'inclusion suivante
   %aaaaaaa
   $$
      Gal(\F) \subset \bigcap_{\I \in spec^{diff}(\O_{IP(\F)})} Stab(\I) 
   $$
   %aaaaaaa
   o\`u $spec^{diff}(\O_{IP(\F)})$ est l'ensemble des id\'eaux diff\'erentiels premiers de $\O_{IP(\F)}$.
   Nous avons d\'ej\`a calcul\'e le $\D$-groupo\"ide de Lie $Gal(\F)$. Les stabilisateurs des id\'eaux diff\'erentiels de $\O_{IP(\F)}$ doivent donc 
   contenir le groupo\"ide de Lie d'invariance de la 2-forme ferm\'ee d\'efinissant le feuilletage.
   
   L'existence de deux int\'egrales premi\`eres ind\'ependantes dans une extension diff\'erentielle du type r\'eductible (d\'efinition \ref{reductibilite})
   assure l'existence d'un id\'eal diff\'erentiel ``assez gros'' dans $\O_{IP(\F)}$. En effet, consid\'erons une telle extension $\CC(x,y,y') \subset K$ 
   et l'application $\O_{IP(\F)} \to K$ qui envoie $H_1,H_2$ sur les deux int\'egrales premi\`eres. Le lemme \ref{extension} nous assure que le type 
   de $K$ est lin\'eaire. Le type de $IP(\F)$ \'etant quadratique, l'application pr\'ec\'edente \`a un noyau non trivial que nous noterons $\I$. 
   Le type de la $\D$-vari\'et\'e d\'efinie par $\I$ est lin\'eaire.

   Consid\'erons le stabilisateur de $\I$, $Stab(\I)$, et le $\D$-groupo\"ide de Lie des transformations de $\CC^3$ fibr\'ees sur l'identit\'e en $x$, 
   $Inv(x)$. La projection sur l'axe des $x$ \'etant transverse au feuilletage, le $\D$-groupo\"ide de Lie $Aut(\F) \cap Inv(x)$ agit simplement transitivement
   sur la $\D$-vari\'et\'e $IP(\F)$. Par cons\'equence le $\D$-groupo\"ide de Lie $Stab(\I)\cap Inv(x)$ agit simplement 
   sur la $\D$-vari\'et\'e d\'efinie par $\I$. Son type est donc lin\'eaire. 
   Le groupo\"ide $Stab(\I) \cap Inv(x)$ est donc strictement plus petit que $Gal(\F) \cap Inv(x)$, ce qui est impossible d'apr\`es la d\'efinition du
   groupo\"ide de Galois. Ceci prouve le th\'eor\`eme. 
\end{preuvedt}

\appendix

\section{Classification de Cartan ; preuve du th\'eor\`eme \ref{itchiban}}

Dans cette annexe, nous reprenons les arguments de \cite{lessousgroupesdesgroupes}, pp 134 --194, pour prouver le th\'eor\`eme suivant. \\

\noindent
%--------------%
\begin{thmref}{itchiban}
%--------------%
\textit{
%--------------%
    Soit $\F$ un feuilletage de $\CC^n$ de codimension deux, d\'efini par une 2-forme ferm\'ee $\gamma$. Le $\D$-groupo\"ide de Lie $Inv(\gamma)$
    d'invariance de cette forme est un $\D$-groupo\"ide de Lie admissible pour $\F$. Si le groupo\"ide de Galois de $\F$ n'est pas $Inv(\gamma)$ alors
    on est dans un des cas suivant :
    \begin{itemize}
       \item $\F$ admet une int\'egrale premi\`ere rationnelle,
       \item il existe une 1-forme alg\'ebrique int\'egrable s'annulant sur $\F$,
       \item il existe un vecteur de 1-forme $\Omega^0 = \left( \begin{array}{c} \omega^0_1 \\ \omega^0_2 \end{array}\right)$ d\'efinissant le 
       feuilletage et une matrice de forme $\Omega^1$ de trace nulle tels que:
       $$
       d\Omega^0 = \Omega^1 \wedge \Omega^0 \ \ \textrm{ et } \ \  d\Omega^1 = \Omega^1 \wedge \Omega^1.
       $$
    \end{itemize}
    }
    \end{thmref}
      Soient $x_0$ un point r\'egulier du feuilletage, $\Theta_\F|_{x_0}$ la forme transverse d'ordre $q$ de $\F$ et 
      $\theta_i^{a,b}$ ses composantes dans la base monomiale sur $J_q(\chi(\CC^2,0))$. 
      Les \'equations de structures s'\'ecrivent
      %aaaaaaaa 
      $$
      d \theta_i^{a,b} = \sum_{ a_1+a_2 = a \atop b_1+b_2=b}{a \choose a_1} {b \choose b_1} \theta_1^{a_1,b_1} \wedge \theta_2^{a_2+1,b_2}
      + \sum_{a_1+a_2 = a \atop b_1+b_2=b} {a \choose a_1} {b \choose b_1} \theta_2^{a_1,b_1} \wedge \theta_2^{a_2,b_2+1}. 
      $$
      pour $a+b \leq q-1$. Nous compl\'etons ces formes en une base des formes transverses au feuilletage sur $Aut(\F)_{q+1}|_{x_0} = R_{q+1}(\F)$
      avec des formes $\omega_i^{a,b}$ pour $a+b = q+1$ et $i=1,2$ satisfaisant les \'equations de structures d'ordre $q+1$ :
      $$
      d\theta_i^{a,b} =  \omega_i^{a+1,b} \wedge \theta_1^{0,0} + \omega_i^{a,b+1} \wedge \theta_2^{0,0} +  \dots,\ \ a+b = q
      $$
      %aaaaaaa
      Ces formes $\omega$ sont ne bien d\'efinies que modulo les formes d'ordre $0$. Les formes $\theta$ \'etant invariantes, les formes $\omega$ ne le sont que 
      modulo les formes d'ordre $0$. 
      %aaaaaaa
      Soit $\Y$ un $\D$-groupo\"ide de Lie transitif admissible pour $\F$ et $p$ l'entier minimal tel que $\Y_p \not = Inv(\gamma)_p$. Nous allons examiner
      les \'equations de structures satisfaites par les formes invariantes sur certaines orbites de $\Y_p|_{x_0}$ et en d\'eduire successivement que 
      
      \begin{itemize}
        \item $p \leq 2$, 
        \item si $p=2$, les \'equations de structure sont de type $asl_2$,
        \item si $p=1$ il existe une combinaison lin\'eaire des formes transverse d'ordre $0$ int\'egrable en restriction \`a une orbite.
      \end{itemize}
      
      \begin{lemme}
       Quitte \`a modifier la suite des formes de Cartan, les restrictions des formes sur une orbite de $Inv(\gamma)|_{x_0}$ satisfont 
      $\theta_1^{a+1,b} + \theta_2^{a,b+1} = 0$.
      \end{lemme}
      \begin{preuve} 
      La forme $t^*\gamma$ est une forme invariante d'ordre 0. Elle est donc \'egale \`a $h\, \theta_1^{0,0} \wedge \theta_2^{0,0}$ o\`u $h$ est 
      un invariant diff\'erentiel d'ordre 1 de $Inv(\gamma)$  donc constant en restriction \`a une orbite. La forme 
      $\theta_1^{0,0} \wedge \theta_2^{0,0}$ est donc ferm\'e ce qui implique $\theta_1^{1,0} + \theta_2^{0,1} = h \theta_1^{0,0} + k\theta_2^{0,0}$ avec $h$ et $k$ des constantes. 
      En faisant agir le diff\'eomorphisme $f(x_1,x_2) = (x_1+hx_1^2,x_2+kx_2^2)$ par composition \`a la source en $x_0$, on trouve une nouvelle suite de formes invariantes commen\c{c}ant par les m\^emes 
      formes d'ordre $0$ et dont les formes d'ordre $1$ sont 
      $$
      \theta_1^{1,0}-h\theta_1^{0,0}, \theta_1^{0,1} + k \theta_1^{0,0}, \theta_2^{1,0} +h \theta_2^{0,0}, \theta_2^{0,1} -k \theta_2^{0,0}.
      $$
      La nouvelle suite que nous noterons encore $\theta$ v\'erifie donc $\theta_1^{1,0} + \theta_2^{0,1} = 0$
      La restriction \`a $Inv(\gamma)|_{x_0}$ de la diff\'erentielle ext\'erieure de cette forme est aussi nulle. Celle-ci est \'egale \`a 
      $$
      (\theta_1^{2,0} + \theta_2^{1,1})\wedge \theta_1^{0,0} + (\theta_1^{1,1} + \theta_2^{0,2})\wedge \theta_2^{0,0}.
      $$  
      On en d\'eduit que $\theta_1^{2,0} + \theta_2^{1,1} = h \theta_2^{0,0}$ et $\theta_1^{1,1} + \theta_2^{0,2} = -h\theta_1^{0,0}$. 
      En faisant agir le diff\'eo\-mor\-phis\-me $f(x_1,x_2) = (x_1+\frac{h}{2}x_2x_1^2, x_2 - \frac{h}{2}x_1x_2^2)$ en $x_0$, on trouve une nouvelle suite
      de formes invariantes satisfaisant les identit\'es annonc\'ees \`a l'ordre $2$. 
      On construit la suite par r\'ecurrence. Comme 
      $$
      d(\theta_1^{a+1,b} + \theta_2^{a,b+1}) = (\theta_1^{a+2,b} + \theta_2^{a+1,b+1})\wedge \theta_1^{0,0} + (\theta_1^{a+1,b+1} + \theta_2^{a,b+2})\wedge \theta_2^{0,0} + \ldots,
      $$ 
      si le premier membre est nul, on a $\theta_1^{a+2,b} + \theta_2^{a+1,b+1} = h \theta_2^{0,0}$ et $\theta_1^{a+1,b+1} + \theta_2^{a,b+2} = -h\theta_1^{0,0}$.
      l'action d'un diff\'eomorphisme en $x_0$ tangent \`a l'ordre $a+b+2$ \`a l'identit\'e permet de construire une suite v\'erifiant les identit\'es anonc\'ees \`a l'ordre $a+b+2$.
      \end{preuve}
      
      Soit maintenant $\Y$ un sous-$\D$-groupo\"ide de Lie transitif de $Inv(\gamma)$ admissible pour $\F$.  
      %
      %-----------------%
      \begin{lemme}
      %-----------------%
      \label{odre<=2}
      %-----------------%
         Si un tel $\D$-groupo\"ide de Lie n'admet pas d'\'equation du premier ordre suppl\'ementaire alors
         il admet quatre \'equations supplementaires d'ordre deux.
      %-----------------%
      \end{lemme}
      %-----------------%
      %
      %----------------%  
      \begin{preuve}
      %----------------%
         Soit $q$ l'ordre minimale des \'equations suppl\'emantaires. Consid\'erons sur $J^*_q|_{x_0}$ la base $\theta_j^\alpha , 
         \ 1\leq j \leq n \ 0 \leq |\alpha| \leq q-1$ des formes (d'ordre $q-1$) invariantes et 
         $\omega_j^\alpha , \ 1\leq j \leq n \ |\alpha| = q$ des formes compl\'etant les $\theta$ en une base de formes transverses. 
         La condition $\theta_1^{a+1,b} + \theta_2^{a,b+1} = 0$ sur $Inv(\gamma)|_{x_0}$ implique que l'on peut choisir des formes $\omega$ 
         satisfaisant $\omega_1^{a+1,b} + \omega_2^{a,b+1} = 0$ pour $a+b=q-1$. 

         Soit $E$ un invariant diff\'erentiel d'ordre $q>1$. La forme $dE = B\omega_2^{q,0} + \sum_{i=0}^q A_i \omega_1^{i,q-i} + ...$ (en 
         n'\'ecrivant pas les formes d'ordre $q-1$) est invariante sous l'action de $\Y_{q+1}$. Les formes 
         $\omega$ \'etant invariantes modulo $\theta_1^{0,0}$ et $\theta_2^{0,0}$, les coefficients $B$ et $A_i$ sont des invariants diff\'erentiels de $\Y$.
         Nous allons maintenant calculer certains termes de la forme $ddE$ restreint sur $\Y_q|_{x_0}$ \`a partir des
         \'equations de structures.
         Dans un premier temps, regardons les termes contenant des formes d'ordre $1$.  Ils sont de la forme 
         $\gamma_1 \wedge \theta_1^{1,0} + \gamma_2 \wedge \theta_1^{1,0}  + \gamma_3 \wedge \theta_2^{1,0}$. Les formes $\gamma_i$
         s'expriment en fonctions des formes $\theta$ et $\omega$. Les termes d'ordre $q$ de ces formes proviennent de la 
         diff\'erentielle des termes d'ordre $q$ de $dE$ que l'on calcule grace aux \'equations de structures. On obtient :
         %aaaaaaaaa
         $$
         \begin{array}{rcrcrcrcrl}
            \gamma_1 & = & (q+1)B \omega_2^{q,0}  &+& (q-1)A_q \omega_1^{q,0}  &+& \sum (2i-q-1)A_i\omega_1^{i,q-i}      &-& (q+1)A_0 \omega_1^{0,q} &+ \ldots \\
            \gamma_2 & =  &  -A_q \omega_2^{q,0}   &+& 2 A_{q-1}\omega_1^{q,0} &+& \sum (q-i+2) A_{i-1}\omega_1^{i,q-i}&+& (q+1)A_0 \omega_1^{1,q-1} &+ \ldots\\ 
            \gamma_3 & = &                                    &-& (q+1)B \omega_1^{q,0}     &+& \sum (i+1) A_{i}\omega_1^{i,q-i}      &+& A_1 \omega_1^{0,q} &+ \ldots
         \end{array}.
         $$
         %aaaaaaaaa   
         Comme $ddE=0$, on en d\'eduit l'existence de trois formes identiques au formes $\gamma_i$ modulo les formes d'ordre $0$ et $1$ s'annulant sur ${\Y_q}|_{x_0}$.
          Comme elles sont in\-d\'e\-pen\-dan\-tes de $dE$, elles 
         proviennent d'invariants diff\'erentiels suppl\'ementaires. En refaisant le m\^eme calcul pour chacune de ces nouvelles formes, 
         on obtient qu'il existe une forme s'annullant sur ${\Y_q}|_{x_0}$ ayant un $A_i$ non nul. On en d\'eduit ensuite que l'on peut supposer que ce coefficient est $A_0$ puis 
         qu'il existe $q+2$ formes s'annullant sur $\Y_q|_{x_0}$ donc $q+2$ invariants diff\'erentiels d'ordre $q$ ($q+3$ invariants d'ordre $q$ donnent un invariant
         d'ordre $q-1$).
         Montrons maintenant que $q$ est inf\'erieur ou \'egale \`a deux. Consid\'erons les $q+2$ invariants diff\'erentiels d'ordre $q$, $E_i$, et leurs 
         diff\'erentielles $dE_i = \sum A_i^j \omega_1^{j,q-j} + B_i \omega_2^{q,0} + \ldots $ o\`u les $A_i^j$ sont des invariants diff\'erentiels de $\Y$.
         Choisissons les combinaisons lin\'eaires de ces formes :
         $$
         \gamma_i = \omega_1^{i,q-i} + \ldots  \text{ et } \gamma_{-1} = \omega_2^{q,0} + \ldots
         $$
          Elles forment un syst\`eme compl\`etement int\'egrable de formes s'annulant sur les orbites de $\Y$ et sont de plus invariantes sous l'action de $\Y$.
          Ces deux conditions impliquent que $d\gamma_i = \sum c_i^{j,k}\gamma_j \wedge\gamma_k$. En r\'eduisant cette \'egalit\'e modulo les formes d'ordre inf\'erieur
          ou \'egal \`a $q-1$, on obtient $0=\sum c_i^{j,k}\omega_1^{j,q-j} \wedge\omega_1^{k,q-k} + \sum c_i^{j,-1}\omega_1^{j,q-j} \wedge\omega_2^{q,0}$. 
          Ces formes \'etant ind\'ependantes, les coefficients $c$ sont tous nuls. Les formes $\gamma$ sont donc ferm\'ees. 
         Regardons maintenant la diff\'erentielle de $\gamma_q$. Cette forme contient le terme 
         $\theta_2^{2,0} \wedge \frac{(q-1)(q-2)}{2}\theta_1^{q-2,1}$. La deuxi\`eme forme de ce produit ext\'erieur ne se 
         repr\'esente dans $d \gamma_q$ que multipil\'ee par une forme d'ordre $1$. La forme $\gamma_q$ \'etant ferm\'ee, ce produit ext\'erieur doit \^etre nul.
         Ceci implique que $q$ est inf\'erieur ou \'egale \`a deux. Si il existe un invariant diff\'erentiel il est donc 
         d'ordre inf\'erieur \`a deux et s'il est d'ordre deux, il y en a quatre.     
      %----------%
      \end{preuve}
      %----------%
      %
      %
      %-----------%
      \begin{lemme}
      %-----------%
         S'il existe un $\D$-groupo\"ide de Lie admissible pour $\F$ ayant un invariant d\'ordre deux
        alors le feuilletage est d\'efini par un vecteur de formes $\Omega^0$ et il existe une matrice de formes $\Omega^1$ de 
         trace nulle telle que $d\Omega^0 = \Omega^1 \wedge \Omega^0$ et $d\Omega^1 = \Omega^1 \wedge \Omega^1$.
     %----------%
      \end{lemme}
     %----------%
     %
     %
     %-----------%
     \begin{preuve}
     %-----------%
         D'apr\'es le lemme pr\'ec\'edent, il y a quatre invariants diff\'erentiels d'ordre 2. Consid\'erons les formes donn\'ees par le lemme pr\'ec\'edent et d\'ecrivant $\Y_2|_{x_0}$.
         %aaaaaaaaa
         $$
         \begin{array}{rl}
         \gamma_1 = & \omega_2^{2,0} + a_1 \theta_2^{1,0} + b_1 (\theta_1^{1,0} - \theta_2^{0,1}) + c_1 \theta_1^{0,1}  + \ldots\\
         \gamma_2 = & \omega_1^{2,0} + a_2 \theta_2^{1,0} + b_2 (\theta_1^{1,0} - \theta_2^{0,1}) + c_2 \theta_1^{0,1} + \ldots \\
         \gamma_3 = & \omega_1^{1,1} + a_3 \theta_2^{1,0} + b_3 (\theta_1^{1,0} - \theta_2^{0,1}) + c_3 \theta_2^{1,0} + \ldots \\
         \gamma_4 = & \omega_1^{0,2} + a_4 \theta_2^{1,0} + b_4 (\theta_1^{1,0} - \theta_2^{0,1}) + c_4 \theta_2^{1,0} + \ldots
         \end{array}
         $$
         %aaaaaaaaa
         Ces formes \'etant invariantes, les coefficients $a$, $b$, $c$ sont des invariants diff\'erentiels de $\Y$.
         Ces formes formant un syst\`eme compl\`etement int\'egrable, elles doivent \^etre ferm\'ees comme nous l'avons montrer au cours de la preuve du lemme pr\'ec\'edent.
         En calculant $d\gamma_1$ modulo les formes d'ordre $0$, on obtient
         %aaaaaaaaa
         $$
         \begin{array}{rl}
         d\gamma_1 = & \frac{3}{2}(\theta_1^{1,0}-\theta_2^{0,1}) \wedge \omega_2^{2,0} - 3\, \theta_2^{1,0}\wedge \omega_1^{2,0} 
                             + a_1 (\theta_1^{1,0} - \theta_2^{0,1}) \wedge \theta_2^{1,0} + da_1 \wedge \theta_2^{1,0} \\
                            & + 2b_1 \theta_2^{1,0}\wedge \theta_1^{0,1} + db_1 \wedge (\theta_1^{1,0} - \theta_2^{0,1}) 
                             + c_1 \theta_1^{0,1} \wedge (\theta_1^{1,0} - \theta_2^{0,1})+dc_1 \wedge \theta_1^{0,1}
         \end{array}
         $$
         %aaaaaaaa
         En examinant les termes contenant des formes d'ordre $2$, on obtient les \'egalit\'es $d a_1 = -3 \gamma_2$, $d b_1 = \frac{3}{2} \gamma_1$ et $dc_1=0$. 
         En faisant le remplacement dans l'\'egalit\'e pr\'ec\'edente, on obtient les \'egalit\'es $c_1 = 0$, $b_2 = \frac{-1}{6}a_1$
         et $c_2 = \frac{-2}{3} b_1$. Les m\^emes manipulations sur
         \begin{itemize} 
         \item $d\gamma_2$ donnent $da_2 = 2\gamma_3$, $db_2 = \frac{1}{2} \gamma_2$ et $dc_2 = -\gamma_1$ 
                puis $b_3 = \frac{-1}{4} a_2$ et $c_3=\frac{-2}{3}a_1$,   
         \item $d\gamma_3$ donnent $da_3 = \gamma_4$, $db_3 = \frac{-1}{2}\gamma_2 $ et $dc_3 = 2 \gamma_2$
                puis $a_4=0$, $b_4 = \frac{-3}{2}a_3$ et $c_4= \frac{3}{2} a_2$,
         \end{itemize}
         On en d\'eduit que, restreint sur la sous-vari\'et\'e invariante $Z_2 =\{ a_1=b_1=a_2=a_3=0\}$, les coefficients $a$, $b$ et $c$ sont nuls.
         En faisant un bon choix des formes $\omega$, les formes restreintes s'\'ecrivent alors :             
         %aaaaaaaaa
         $$
         \begin{array}{rl}
         \gamma_1 = & \omega_2^{2,0} \\
         \gamma_2 = & \omega_1^{2,0} + h_2 \theta_1^{0,0}\\
         \gamma_3 = & \omega_1^{1,1} + h_3 \theta_1^{0,0} \\
         \gamma_4 = & \omega_1^{0,2} + h_4 \theta_1^{0,0}.
         \end{array}
         $$
         %aaaaaaaaa
         Les formes $\theta$ restreintes satisfont les \'equations :
         $$
         \begin{array}{rl}
         d(\theta_1^{1,0}-\theta_2^{0,1}) & =  2 \theta_2^{1,0} \wedge \theta_1^{0,1} + (h_2-h_1) \theta_2 \wedge \theta_1\\
         d\theta_1^{0,1} & =   \theta_1^{0,1} \wedge (\theta_1^{1,0}-\theta_2^{0,1}) + h_3 \theta_2 \wedge \theta_1\\
         d\theta_2^{1,0} & = - \theta_2^{1,0} \wedge (\theta_1^{1,0}-\theta_2^{0,1}) - h_1 \theta_2 \wedge \theta_1 .\\
         \end{array}
         $$
         En v\'erifiant que ces formes sont ferm\'ees, on obtient $h_1=h_2=h_3=0$.
         Les formes $\theta$ restreintes sur $Z_2 = \{ a_1=b_1=a_2=a_3=0\}$ v\'erifient les identit\'es annonc\'ees dans la proposition.
         En choisissant ensuite une section $f$ de $Z_2$ pour la projection but, on tire les formes sur $\CC^n$ : $\Omega = f^*\Theta$. 
      %-----------------%
      \end{preuve}
      %-----------------%

      %
      %-----------------%
      \begin{lemme}
      %-----------------%
      %
      \label{lemmedeitchiban}
      %-----------------%
         Si un $\D$-groupo\"ide de Lie admissible pour $\F$ admet des \'equations du premier ordre suppl\'ementaires alors il laisse un feuilletage de 
         codimension un invariant.  
      %-----------------%
      \end{lemme}
      %-----------------%
      %
      %-----------------%
      \begin{preuve}
      %-----------------% 
      Nous allons examiner $\Y_1|_{x_0}$ muni des deux formes invariantes $\theta_1^{0,0}$ et $\theta_2^{0,0}$ que l'on compl\`ete
      en une base avec des formes $\omega_j^{k,\ell}$ pour $j=1,2$ et $k+\ell=1$ satisfaisant les premi\`eres \'equations de structure.
      
      Supposons que $\Y_1$ est d\'efinie par un invariant d'ordre un suppl\'ementaire $E$. On \'ecrit 
      %aaaaaaaa
      $$
      dE = a \omega_2^{1,0}+b (\omega_1^{1,0} - \omega_2^{0,1}) + c \omega_1^{0,1}  + h \theta_1^{0,0} + k \theta_2^{0,0}
      $$ 
      %aaaaaaaa
      En choisissant correctement les formes $\omega$, on peut supposer $k_1=k_2=0$. Les coefficients $a$, $b$ et $c$ \'etant des invariants, 
      la forme $\gamma = \omega_2^{1,0}+b (\omega_1^{1,0} - \omega_2^{0,1}) + c \omega_1^{0,1} $, obtenue en divisant par $a$, 
      est ferm\'ee. En calculant $d\gamma$ modulo $\gamma$, on obtient $b^2+c =0$. 
      En \'ecrivant ensuite que $\gamma$ est ferm\'ee, on obtient 
      %aaaaaaaa
      $$ 
      0= db +\omega_2^{1,0}+b (\omega_1^{1,0} - \omega_2^{0,1}) - b^2 \omega_1^{0,1} .
      $$
      %aaaaaaaa 
      En restriction \`a la sous-vari\'et\'e $Z_1=\{b=0\}$, la forme $\theta_2^{0,0}$ est int\'egrable. Le pull-back de cette forme par 
      une section de $Z_1$ donne une forme alg\'ebrique int\'egrable s'annulant sur le feuilletage.  
       
      Supposons que $\Y_1$ soit donn\'e par deux invariants $E_1$ et $E_2$. Les coefficients des formes $\omega$ dans les diff\'erentielles $dE_1$ 
      et $dE_2$ sont des invariants. En prenant des combinaisons lin\'eaires de ces formes, on trouve les deux formes invariantes
      %aaaaaaa
      $$
      \begin{array}{rl}
      \gamma_1 & = \omega_2^{1,0} + a\, \omega_1^{0,1} + h_1 \theta_1^{0,0} + k_1 \theta_2^{0,0} \\
      \gamma_2 & = (\omega_1^{1,0}-\omega_2^{0,1}) + b\, \omega_1^{0,1} + h_2 \theta_1^{0,0} + k_2 \theta_2^{0,0}
      \end{array}
      $$
      %aaaaaaa      
      En choisissant correctement les formes $\omega$, on peut supposer que $h_1=k_1=h_2=k_2=0$.
      Ces formes s'annulant sur les orbites de $\Y$, on a $d\gamma_1 = c_1 \gamma_1 \wedge \gamma_2$ et 
      $d\gamma_2 = c_2 \gamma_1 \wedge \gamma_2$. Le calcul direct \`a partir des formules donne 
      $$
      \begin{array}{rcl}
      d\gamma_1  & = &-\gamma_1 \wedge \gamma_2 + (da -2a \gamma_2 +b \gamma_1) \wedge \omega_1^{0,1} \\
      d\gamma_2  & = & (db + 2\gamma_1 - b \gamma_2) \wedge  \omega_1^{0,1}.
      \end{array}
      $$
      En se pla\c{c}ant sur l'hypersurface $Z_1=\{b=0\}$, la forme $\theta_2^{0,0}$ est int\'egrable.
      
      Supposons enfin que $\Y_1$ soit donn\'ee par trois invariants. En effectuant les manipulations pr\'ec\'edentes, on obtient trois formes invariantes s'annulant sur $\Y_1$
      %aaaaaaa
      $$
      \begin{array}{rcl}
      \gamma_1 & = & \omega_2^{1,0} + h_1 \theta_1^{0,0} + k_1 \theta_2^{0,0} \\
      \gamma_2 & = & (\omega_1^{1,0}-\omega_2^{0,1}) + h_2 \theta_1^{0,0} + k_2 \theta_2^{0,0}\\
      \gamma_3 & = & \omega_1^{0,1} + h_3 \theta_1^{0,0} + k_3 \theta_2^{0,0}
      \end{array}
      $$
      Quitte \`a choisir d'autres formes $\omega$, on peut supposer $h_1=k_1=k_2=k_3=0$. En calculant les diff\'erentiels et en v\'erifiant que l'on
      doit obtenir des formules $d\gamma_i = \sum c_i^{j,k} \gamma_j \wedge \gamma_k$, on obtient $h_2=h_3=0$. En restriction \`a $Y_1$,
      la forme $\theta_2^{0,0}$ est int\'egrable (m\^eme ferm\'ee).  
      %aaaaaaa      
      %-----------------%
      \end{preuve}
      %-----------------%
    Dans le cas o\`u $\Y$ admet un invariant d'ordre 0, le th\'eor\`eme \ref{theoinvdiff} prouve que le feuilletage admet une int\'egrale premi\`ere.

\end{document}